\documentclass{article}
\usepackage{amsfonts}
\usepackage{graphicx}
\usepackage{epstopdf}
\usepackage{color,url,bm}
\usepackage{marginnote}
\usepackage{hyperref}

\newtheorem{theo}{Theorem}
\newtheorem{lem}{Lemma}
\newtheorem{cor}{Corollary}
\newtheorem{rem}{Remark}

\def\pmatrix{\left(\begin{array}}
\def\endpmatrix{\end{array}\right)}
\def\no{\noindent}

\def\proof{\no\underline{Proof}\quad}
\def\QED{\fbox{}}
\def\dd{\mathrm{d}}
\def\diag{\mathrm{diag}}

\def\RR{\mathbb{R}}

\def\P{{\cal P}}
\def\I{{\cal I}}
\def\M{{\cal M}}

\def\bfe{{\bm{e}}}
\def\bfg{\bm{g}}
\def\bfv{{\bm{v}}}
\def\bfpsi{{\bm{\psi}}}
\def\bflam{{\bm{\lambda}}}

\def\bfzero{{\bm{0}}}
\def\bfgamma{{\bm{\gamma}}}
\def\bfzeta{{\bm{\zeta}}}

\def\hH{\hat{H}}

\title{Arbitrarily high-order energy-conserving methods for Hamiltonian problems with quadratic holonomic constraints}
\author{P.\,Amodio\,\thanks{Dipartimento di Matematica, Universit\`a di Bari, Italy.} \and L.\,Brugnano\,\thanks{Dipartimento di Matematica e Informatica ``U.\,Dini'', Universit\`a di Firenze, Italy.} \and G.\,Frasca-Caccia\,\thanks{Dipartimento di Matematica, Universit\`a di Salerno, Italy.} \and F.\,Iavernaro\,\footnotemark[1]}

\begin{document}
\maketitle

\begin{abstract} In this paper, we define arbitrarily high-order energy-conserving methods for Hamiltonian systems with quadratic holonomic constraints.  The derivation of the methods is made within the so-called {\em line integral framework}. Numerical tests to illustrate the theoretical findings are presented.

\smallskip
\no{\bf Keywords:} constrained Hamiltonian systems; quadratic holonomic constraints; energy-conserving methods; line integral methods; Hamiltonian Boundary Value Methods; HBVMs.\\

\smallskip
\no{\bf MSC:} 65P10, 65L80, 65L06.
\end{abstract}

\section{Introduction}\label{intro}
In recent years, much interest has been given to the modeling and/or simulation of tethered systems, where the dynamics of interconnected bodies is studied (see, e.g. \cite{Kaw2011,Xu2012,Jung2014,Jung2015,Zab2015,Zakr2015,Wang2016,Yu2016,Baiao2017,Lian2019,Luo2021}). It turns out that the underlying dynamics is often described by a Hamiltonian system, for which the total energy is conserved.

Motivated by this fact, we here investigate the numerical approximation of a constrained Hamiltonian dynamics, described by the separable Hamiltonian
\begin{equation}\label{Hqp} 
H(q,p) =  \frac{1}2 p^\top  M^{-1} p -U(q), \qquad q,p\in\RR^m,
\end{equation}
where $M$ is a symmetric and positive-definite (spd) matrix, subject to $\nu$ quadratic {\em holonomic} constraints,
\begin{equation}\label{gq}
g(q) = 0\in\RR^\nu, \qquad \nu\le m,
\end{equation}
i.e., the entries of $g$  are quadratic polynomials. Hereafter, we shall assume all points be regular for the constraints, i.e., $\nabla g(q)\in\RR^{m\times \nu}$ has full column rank or, equivalently, 
\begin{equation}\label{regular}
\nabla g(q)^\top M^{-1}\nabla g(q) \in\RR^{\nu\times\nu}  \qquad \mbox{is~spd}.
\end{equation}
Moreover, we shall assume that its smallest eigenvalue is bounded away from 0, in the domain of interest. Also,   
for sake of simplicity,  in the same domain the potential $U$ will be assumed to be analytic.

 It is well-known that the problem defined by (\ref{Hqp})--(\ref{gq}) can be cast in Hamiltonian form by defining the augmented Hamiltonian
\begin{equation}\label{Hqpl}
\hH(q,p,\lambda) = H(q,p) + \lambda^\top g(q),
\end{equation}
where $\lambda$ is the vector of the Lagrange multipliers. The resulting constrained Hamiltonian system reads: 
\begin{equation}\label{constrp}
\dot q =  M^{-1}p, \qquad \dot p = \nabla U(q)-\nabla g(q)\lambda, \qquad g(q)=0, \qquad  t\in[0,T],
\end{equation}
and is subject to consistent initial conditions,
\begin{equation}\label{q0p0}
q(0)=q_0,\qquad p(0)=p_0,
\end{equation}
such that
\begin{equation}\label{consist0}
g(q_0)=0, \qquad \nabla g(q_0)^\top  M^{-1} p_0 = 0.
\end{equation}
Clearly, $H(q,p)\equiv\hH(q,p,\lambda)$, provided that the constraints (\ref{gq}) are satisfied, and a straightforward calculation proves that both are conserved along the solution trajectory.

We notice that the condition $g(q_0)=0$ ensures that $q_0$ belongs to the manifold
\begin{equation}\label{M}
\M = \left\{ q\in\RR^m: g(q)=0\right\},
\end{equation}
as required by the constraints, whereas the condition $\nabla g(q_0)^\top  M^{-1}p_0$ means that the motion initially stays on the tangent space to $\M$ at $q_0$. This condition is satisfied by all points on the solution trajectory, since, in order for the constraints to be conserved,  the following condition needs to be satisfied as well:
\begin{equation}\label{Dgp0}
\dot g(q) = \nabla g(q)^\top \dot q = \nabla g(q)^\top  M^{-1} p = 0 \in\RR^\nu.
\end{equation}
These latter constraints are usually referred to as {\em hidden constraints}, and allow the derivation of the vector of the Lagrange multiplier $\lambda$. In fact, from (\ref{Dgp0}) and (\ref{constrp})-(\ref{q0p0}), one obtains
\begin{eqnarray}\nonumber
0&=&\nabla g(q(t))^\top  M^{-1} p(t) \\
&=& \nabla g(q(t))^\top  M^{-1}\left[ p_0 + \int_0^t \nabla U(q(\zeta))\dd\zeta - \int_0^t \nabla g(q(\zeta))\lambda(\zeta)\dd\zeta\right],\label{hiddenp}
\end{eqnarray}
from which one derives the integral equation:
\begin{eqnarray}\nonumber\lefteqn{
\nabla g(q(t))^\top  M^{-1}\int_0^t \nabla g(q(\zeta))\lambda(\zeta)\dd\zeta} \\
&=& \nabla g(q(t))^\top  M^{-1} \left[ p_0 + \int_0^t \nabla U(q(\zeta))\dd\zeta\right]. \label{inteq}
\end{eqnarray}

\begin{rem}\label{lam}
From (\ref{hiddenp})-(\ref{inteq}), one deduces that the vector of the Lagrange multipliers depends on the functions $q$ and $p$. We shall denote this by using the notation:
\begin{equation}\label{lamq}
\lambda(t) := \lambda(q(t),p(t)).
\end{equation}
\end{rem}

We stress that the condition (\ref{Dgp0}) can be conveniently relaxed for a numerical approximation: in fact, we only ask for $\nabla g(q)^\top  M^{-1}p$ to be suitably small along the numerical solution. Consequently, when solving the problem on the interval $[0,h]$, we require that the approximations,
\begin{equation}\label{q1p1}
q_1\approx q(h), \qquad p_1\approx p(h),
\end{equation} 
satisfy the conservation of both the Hamiltonian and the constraints, 
\begin{equation}\label{Hgcons}
H(q_1,p_1) = H(q_0,p_0),  \qquad g(q_1)=0,
\end{equation}
whereas the hidden constraints, see (\ref{Dgp0}), become:
\begin{equation}\label{hidcons}
\nabla g(q(h))^\top  M^{-1}p(h) = O(h^{\bar r}),
\end{equation}
for a convenient $\bar r\ge1$.

For sake of completeness, we recall that a formal expression for the vector $\lambda$ is obtained by further differentiating (\ref{Dgp0}):
\begin{equation}\label{lam0}
\ddot g(q) = \nabla^2 g(q) (M^{-1}p,M^{-1}p) + \nabla g(q)^\top M^{-1} \left[ \nabla U(q) - \nabla g(q)\lambda\right].
\end{equation}
In fact, by imposing the vanishing of this derivative yields 
\begin{equation}\label{lambda1}
\left[\nabla g(q)^\top M^{-1} \nabla g(q)\right]\lambda ~=~  \nabla^2 g(q) (M^{-1}p,M^{-1}p) + \nabla g(q)^\top M^{-1} \nabla U(q).
\end{equation}
Consequently, the following result follows. 

\begin{theo}\label{lamex} The vector $\lambda$ exists and is uniquely determined, provided that (\ref{regular}) holds true. 
In fact, in such a case, from (\ref{lambda1}), we obtain
$$
\lambda =  \left[\nabla g(q)^\top M^{-1} \nabla g(q)\right]^{-1}\left[ \nabla^2 g(q) (M^{-1}p,M^{-1}p) + \nabla g(q)^\top M^{-1} \nabla U(q)\right].
$$
\end{theo}


The numerical solution of Hamiltonian problems with holonomic constraints has been the subject of several investigations across the years. Among the proposed approaches, we mention: the popular Shake-Rattle method \cite{Shake1977,Rattle1983}, later proved to be symplectic \cite{LeSk1994}; higher order methods obtained via symplectic PRK methods \cite{Jay1996}; composition methods \cite{Reich1996,Reich1997}; symmetric LMFs \cite{CHL2013}; methods based on discrete derivatives \cite{Gonzalez1999}; local parametrizations of the manifold containing the solution \cite{BCF2001}; projection techniques \cite{Seiler1998,WDLW2016}. As additional references, we refer to \cite{LeRe1994,Seiler1999,Hairer2003,LBS2004,LW2022} and to the monographs \cite{BCP1996,LeRe2004,HLW2006,HaWa2010}, and references therein.

In this paper we pursue a different approach, utilizing the so-called {\em line integral framework}, which has already been used for deriving the energy-conserving Runge-Kutta methods, for unconstrained Hamiltonian systems, named {\em Hamiltonian Boundary Value Methods (HBVMs)} \cite{BIT2009,BIT2010,BIT2011,BIT2012,BIT2015}, see also the  monograph \cite{LIMbook2016} and the review paper \cite{axioms2018}. This approach has been already used in \cite{BGIW2018} for deriving energy-conserving methods for constrained Hamiltonian problems which, however, are in general only second-order accurate. In this paper we propose a class of arbitrarily high-order energy conserving methods for constrained Hamiltonian dynamics, in the relevant case where the constraints are quadratic. As matter of fact, the proposed methods provide approximations (\ref{q1p1}) satisfying (\ref{Hgcons}) and(\ref{hidcons}) with $\bar r$ arbitrarily large.

With this premise, the structure of the paper is as follows: in Section~\ref{polapp} we define the polynomial approximation defining the basic step of application of the method also studying its conservation properties; in Section~\ref{order} we study the order of convergence of the method; in Section~\ref{discrete} we discuss its fully discrete implementation; numerical examples are reported in Section~\ref{numex}; at last, some concluding remarks are given in Section~\ref{fine}.

\section{Polynomial approximation}\label{polapp}
Let us consider the following orthonormal polynomial basis,\footnote{As is usual, $\Pi_j$ is the vector space of polynomials of degree $j$, and $\delta_{ij}$ is the Kronecker symbol.}
\begin{equation}\label{leg}
P_j\in\Pi_j, \qquad \int_0^1 P_i(x)P_j(x)\dd x=\delta_{ij}, \qquad \forall i,j=0,1,\dots,
\end{equation}
given by the shifted and scaled Legendre polynomials, and the infinite expansions on the interval $[0,h]$ (see (\ref{constrp})):
\begin{equation}\label{exp}
p(ch) = \sum_{j\ge0} P_j(c)\gamma_j(p), \quad \nabla U(q(ch)) = \sum_{j\ge0} P_j(c)\psi_j(q),\quad c\in[0,1],
\end{equation}
where, for all $j=0,1,\dots$,
\begin{equation}\label{gampsij}
\gamma_j(p) = \int_0^1 P_j(\zeta)p(\zeta h)\dd\zeta, \qquad \psi_j(q) = \int_0^1 P_j(\zeta)\nabla U(q(\zeta h))\dd\zeta.
\end{equation}
In order to derive polynomial approximations $u$ and $v$ of degree $s$,
$$u(ch)\approx q(ch), \qquad v(ch)\approx p(ch), \qquad c\in[0,1],$$
we shall consider the following approximation to (\ref{constrp}) on the interval $[0,h]$,
\begin{eqnarray}\nonumber
\dot u(ch) &=& M^{-1}\sum_{j=0}^{s-1} P_j(c)\gamma_j(v), \qquad u(0)=q_0, \qquad c\in[0,1],\\ \label{lim}
\dot v(ch) &=& \sum_{j=0}^{s-1} P_j(c)\psi_j(u)-\sum_{j=1}^s\ell_j(c) \rho_j, \qquad v(0)=p_0,\\
\rho_j &:=& \nabla g(u(c_jh))\lambda_j, \qquad j=1,\dots,s,\nonumber 
\end{eqnarray}
where:
\begin{itemize}
\item $\gamma_j(v)$ and $\psi_j(u)$ are defined according to (\ref{gampsij}), by formally replacing $q$ and $p$ with $u$ and $v$, respectively;
\item $0< c_1<\dots<c_s< 1$~ are the zeros of ~$P_s(c)$, ~$P_s(c_j)=0$, ~$j=1,\dots,s$;
\item $\ell_j(c)$ are the Lagrange polynomials defined over the previous abscissae;
\item finally, by using the notation (\ref{lamq}),
\begin{equation}\label{blam} 
\lambda_j\approx \bar\lambda(c_jh) := \lambda(u(c_jh),v(c_jh)) 
\end{equation}
are approximation to the Lagrange multipliers, which will be defined in the sequel.
\end{itemize}
Moreover, let us denote by 
$$b_j = \int_0^1 \ell_j(c)\dd c, \qquad j=1,\dots,s,$$
the weights of the Gauss-Legendre quadrature of order $2s$. 
The new approximations,
\begin{equation}\label{q1p1def}
q_1:=u(h), \qquad p_1:=v(h),
\end{equation}
will be then defined such that (\ref{Hgcons}) and (\ref{hidcons}) are satisfied, for a convenient $\bar r$ depending on $s$.

For later use, we observe that, at the Gauss-Legendre abscissae $\{c_i\}$, (\ref{lim}) reads:
\begin{eqnarray}
\dot u(c_ih) &=& M^{-1}v(c_ih),  \label{limci}\\
\dot v(c_ih) &=& \sum_{j=0}^{s-1} P_j(c_i)\psi_j(u)-\nabla g(u(c_ih))\lambda_i, \qquad i=1,\dots,s.\nonumber
\end{eqnarray}
In fact, since $v(ch)\in\Pi_s$, then (see (\ref{gampsij}))
$$v(ch) = \sum_{j=0}^s P_j(c)\gamma_j(v), \qquad c\in[0,1],$$
however, at the zeros of $P_s$, one obviously has:
\begin{equation}\label{vciheq}
v(c_ih) =  \sum_{j=0}^{s-1} P_j(c_i)\gamma_j(v), \qquad i=1,\dots,s.
\end{equation}

\subsection{Constraints conservation}

Let us now study the conservation of the holonomic constraints in (\ref{Hgcons}), i.e.
under which conditions $g(q_1)=0$. One has, by virtue of (\ref{limci}):
\begin{eqnarray*}\nonumber
g(q_1) &\equiv& g(u(h)) ~=~ h\int_0^1 \nabla g(u(ch))^\top \dot u(ch)\dd c \\ 
&=&h \sum_{i=1}^s b_i  \nabla g(u(c_ih))^\top \dot u(c_ih)  ~=~ h \sum_{i=1}^s b_i  \nabla g(u(c_ih))^\top M^{-1} v(c_ih),\qquad\quad
\end{eqnarray*}
where the third equality follows from the fact that the integrand has degree $2s-1$, and the last equality follows from (\ref{limci}).
Consequently, the following result is proved.

\begin{theo}\label{gheq0} $g(q_1)=0$, provided that:
\begin{equation}\label{azzera0}
\nabla g(u(c_ih))^\top M^{-1} v(c_ih) = 0\in\RR^\nu, \qquad i=1,\dots,s.
\end{equation}
\end{theo}




\subsection{Energy conservation} 

Let us now study the conservation of the Hamiltonian function (\ref{Hqp}). By using the usual line integral argument, from (\ref{lim}) one has:

\begin{eqnarray*}\nonumber
\lefteqn{H(q_1,p_1)-H(q_0,p_0)~\equiv~ H(u(h),v(h))-H(u(0),v(0))}\\ \nonumber
&=& \int_0^h \frac{\dd}{\dd t} H(u(t),v(t))\dd t\\ \nonumber
&=& h\int_0^1 \left[H_q(u(ch),v(ch))^\top \dot u(ch)+ H_p(u(ch),v(ch))^\top \dot v(ch)\right]\dd c\\ \nonumber
&=& h\int_0^1\left[-\nabla U(u(ch))^\top M^{-1}\sum_{j=0}^{s-1}P_j(c)\gamma_j(v)+v(ch)^\top M^{-1}\sum_{j=0}^{s-1}P_j(c)\psi_j(u)\right.\\ \nonumber
&&\left.\qquad ~~-~v(ch)^\top M^{-1}\sum_{j=1}^s\ell_j(c)\rho_j\right]\dd c\\ \nonumber
&=& -~h \sum_{j=0}^{s-1} \overbrace{\left[\int_0^1P_j(c)\nabla U(u(ch))\dd c\right]^\top}^{=\psi_j(u)^\top} M^{-1}\gamma_j(v)\\ \nonumber
&& +~h \sum_{j=0}^{s-1} \overbrace{\left[\int_0^1P_j(c)v(ch)\dd c\right]^\top}^{=\gamma_j(v)^\top} M^{-1}\psi_j(u) 
-h\int_0^1 \sum_{j=1}^s v(ch)^\top M^{-1} \ell_j(c)\rho_j\dd c \\ \nonumber
&=& -h\sum_{j=0}^{s-1} \psi_j(u)^\top M^{-1}\gamma_j(v)+h\sum_{j=0}^{s-1}\gamma_j(v)^\top M^{-1}\psi_j(u) -h\sum_{i=1}^s b_i v(c_ih)^\top M^{-1}\rho_i\\ 
&=& -h\sum_{i=1}^s b_i v(c_ih)^\top M^{-1} \nabla g(u(c_ih))\lambda_i. 
\end{eqnarray*}
In fact, considering that the Gauss-Legendre formula has order $2s$, 
\begin{eqnarray*}\lefteqn{
\int_0^1 \sum_{j=1}^s v(ch)^\top M^{-1} \ell_j(c)\rho_j\dd c ~=~ \sum_{i,j=1}^s b_i  v(c_ih)^\top M^{-1}\ell_j(c_i)\rho_j} \\
&=& \sum_{i=1}^s b_i v(c_ih)^\top M^{-1}\rho_i ~=~\sum_{i=1}^s b_i v(c_ih)^\top M^{-1} \nabla g(u(c_ih))\lambda_i.
\end{eqnarray*}
Consequently, the following result is proved. 

\begin{theo}\label{evai} 
Assume that the conditions (\ref{azzera0}) hold true.  Then, $H(q_1,p_1)=H(q_0,p_0)$, i.e., the energy is conserved.
\end{theo}

\subsection{Computing the approximate Lagrange multipliers}

Before studying the accuracy of the approximations (\ref{q1p1}) and (\ref{q1p1def}), let us derive the algebraic formulation of the conditions (\ref{azzera0}), i.e., the discrete analog of (\ref{hiddenp})
at $t=c_ih$, $i=1,\dots,s$. For this purpose, let us define the following matrices,
\begin{eqnarray}\nonumber
\nabla\bfg &=& \pmatrix{ccc} \nabla g(u(c_1h))& &\\ &\ddots\\ &&\nabla g(u(c_sh))\endpmatrix\in\RR^{sm\times s\nu}\\[2mm] \label{matrici}
\P_s &=&\pmatrix{ccc} P_0(c_1)&\dots&P_{s-1}(c_1)\\ \vdots & &\vdots\\ P_0(c_s)& \dots &P_{s-1}(c_s)\endpmatrix\in\RR^{s\times s},\\[2mm] \nonumber
\I_s &=&\pmatrix{ccc} \int_0^{c_1}P_0(x)\dd x&\dots&\int_0^{c_1}P_{s-1}(x)\dd x\\ \vdots & &\vdots\\ \int_0^{c_s}P_0(x)\dd x& \dots &\int_0^{c_s}P_{s-1}(x)\dd x\endpmatrix\in\RR^{s\times s},\\[2mm] \nonumber
\Omega &=&\diag(b_1,\,\dots,\,b_s)\in\RR^{s\times s},
\end{eqnarray}
and vectors,
\begin{eqnarray}\nonumber
\bfpsi &=& \pmatrix{c} \psi_0(u)\\ \vdots\\ \psi_{s-1}(u)\endpmatrix, ~
\bfv    \,=\, \pmatrix{c} v(c_1h)\\ \vdots\\ v(c_sh)\endpmatrix~\in\RR^{sm},\\ \label{vettori}
\bflam &=& \pmatrix{c}\lambda_1\\ \vdots\\ \lambda_s\endpmatrix,~\bfzero \,=\, \pmatrix{c}0\\ \vdots\\ 0\endpmatrix^\top~\in\RR^{s\nu},\\ \nonumber
\bfe     &=& (1,\,\dots,\,1)^\top\in\RR^s.\ \nonumber
\end{eqnarray}
In so doing, it turns out that (\ref{azzera0}) can be cast as:
$$\nabla \bfg^\top I_s\otimes M^{-1}\, \bfv = \bfzero.$$
Moreover, from (\ref{lim}) one derives:
$$\bfv = \bfe\otimes p_0 +h\I_s\otimes I_m \bfpsi -h\I_s\P_s^\top\Omega\otimes I_m \nabla\bfg\cdot\bflam.$$
Consequently, combining the two last equations, one formally obtains:  
\begin{equation}\label{lam}
\left( \nabla \bfg^\top \, h\I_s\P_s^\top\Omega\otimes M^{-1}\, \nabla\bfg\right)\bflam = 
\nabla\bfg^\top \left(  \bfe\otimes M^{-1}p_0 + h\I_s\otimes M^{-1}\bfpsi\right),
\end{equation}
which, as observed above, represents the discrete counterpart of (\ref{inteq}) at $t=c_ih$, $i=1,\dots,s$. Consequently, the following result is proved.

\begin{cor}\label{lam1}
Conditions (\ref{azzera0}) are equivalent to (\ref{matrici})--(\ref{lam}), i.e.,
\begin{eqnarray}\nonumber\lefteqn{
h \nabla g(u(c_ih))^\top M^{-1}  \sum_{\ell=0}^{s-1} \int_0^{c_i}P_\ell(x)\dd x \sum_{j=1}^s b_j P_\ell(c_j)\nabla g(u(c_jh))\lambda_j}\\ \label{azzera1}
&\equiv& h \nabla g(u(c_ih))^\top M^{-1} \sum_{j=1}^s\overbrace{b_j \sum_{\ell=0}^{s-1} \int_0^{c_i}P_\ell(x)\dd x\, P_\ell(c_j)}^{=\,a_{ij}}\nabla g(u(c_jh))\lambda_j\\
&=& \nabla g(u(c_ih))^\top M^{-1}\left[ p_0 + h\sum_{j=0}^{s-1}\int_0^{c_i}P_j(x)\dd x\,\psi_j(u)\right], \qquad i=1,\dots,s.\nonumber
\end{eqnarray}
\end{cor}

We observe that  $a_{ij}$ in (\ref{azzera1}) is nothing but the $(i,j)th$ entry of the Butcher matrix of the $s$-stage Gauss method (see, e.g., \cite{LIMbook2016,axioms2018}). Moreover, the problem (\ref{lam}) is well-posed, since $\nabla\bfg$ has full column rank, by hypothesis, and
\begin{equation}\label{IsPs}
\I_s = \P_sX_s, \qquad \P_s^\top\Omega = \P_s^{-1},
\end{equation}
with
\begin{equation}\label{Xs}
X_s = \pmatrix{cccc} \xi_0&-\xi_1\\ \xi_1 &0 &\ddots\\ &\ddots &\ddots &-\xi_{s-1}\\ & &\xi_{s-1}&0\endpmatrix,\qquad \xi_i=\left(2\sqrt{|4i^2-1|}\right)^{-1},
\end{equation}
a nonsingular matrix. Let us now study, for later use, the differences (see (\ref{blam})):
\begin{equation}\label{dli}
\delta\lambda_j := \bar\lambda(c_jh)-\lambda_j, \qquad j=1,\dots,s.
\end{equation}
For this purpose, we need the following lemma.

\begin{lem}\label{Ghj}
Let $G:[0,h]\rightarrow V$, with $V$ a vector space, admit a Taylor expansion at 0. Then, 
$$\int_0^1 P_j(c) G(ch)\dd c = O(h^j), \qquad j=0,1,\dots.$$
\end{lem}
\proof See \cite[Lemma\,1]{BIT2012}.\,\QED\bigskip

We also recall that, for any suitably regular function $G$,
\begin{eqnarray}\label{quaderr}
\int_0^1 P_j(c)G(ch)\dd c &=& \sum_{i=1}^s b_i P_j(c_i)G(c_ih) \,+\,\Delta_j(h), \\
&& \Delta_j(h)=O(h^{2s-j}), \qquad j=0,\dots,s-1, \nonumber
\end{eqnarray}
$(c_i,b_i)$ being, as usual, the abscissae and weights of the Gauss-Legendre quadrature of order $2s$. 
We are now able to prove the following result.

\begin{theo}\label{dlamj} Concerning the differences (\ref{dli}),
one has:
$$\delta\lambda_j \equiv \bar\lambda(c_jh)-\lambda_j = O(h^s).$$
\end{theo}
\proof
In fact, from (\ref{inteq}) with $q$ and $p$ replaced by $u$ and $v$, respectively, and evaluated at $t=c_ih$, one obtains:
$$
\nabla g(u(c_ih))^\top M^{-1} v(c_ih) = 0,\qquad i=1,\dots,s,
$$
i.e.,  
\begin{eqnarray}\nonumber
\lefteqn{\nabla g(u(c_ih))^\top M^{-1} h \int_0^{c_i} \nabla g(u(\xi h))\bar\lambda(\xi h)\dd\xi} \\ \label{dl1}
&=& \nabla g(u(c_ih))^\top M^{-1}\left[ p_0 + h\int_0^{c_i} \nabla U(u(\xi h))\dd\xi\right], \qquad i=1,\dots,s.\qquad
\end{eqnarray}
Considering that $\nabla g(u)\in\Pi_s$, and by virtue of Lemma~\ref{Ghj} and (\ref{quaderr}), one has:
\begin{eqnarray}\nonumber
\lefteqn{\int_0^{c_i} \nabla g(u(\xi h))\bar\lambda(\xi h)\dd\xi }\\ \nonumber 
&=& \sum_{\ell=0}^{s-1} \int_0^{c_i} P_\ell(x)\dd x \int_0^1 P_\ell(\tau)\nabla g(u(\tau h))\bar\lambda(\tau h)\dd\tau + O(h^s)\\ \nonumber
&=&\sum_{\ell=0}^{s-1} \int_0^{c_i} P_\ell(x)\dd x \left[ \sum_{j=1}^s b_j P_\ell(c_j)\nabla g(u(c_j h))\bar\lambda(c_jh) + O(h^{2s-\ell})\right] + O(h^s)\\ \label{dl2}
&=&\sum_{\ell=0}^{s-1} \int_0^{c_i} P_\ell(x)\dd x \sum_{j=1}^s b_j P_\ell(c_j)\nabla g(u(c_j h))\bar\lambda(c_jh) + O(h^s).
\end{eqnarray}
Similarly, from (\ref{exp})-(\ref{gampsij}), one has:
\begin{eqnarray}\nonumber
\int_0^{c_i} \nabla U(u(\xi h))\dd\xi &=& \sum_{j\ge0} \int_0^{c_i} P_j(x)\dd x \, \psi_j(u) \\
&=& \sum_{j=0}^{s-1} \int_0^{c_i} P_j(x)\dd x \, \psi_j(u) + O(h^s). \label{dl3}
\end{eqnarray}
From (\ref{dl1})--(\ref{dl3}) one then obtains:
\begin{eqnarray} \label{azzera2}
\lefteqn{
h \nabla g(u(c_ih))^\top M^{-1}  \sum_{\ell=0}^{s-1} \int_0^{c_i}P_\ell(x)\dd x \sum_{j=1}^s b_j P_\ell(c_j)\nabla g(u(c_jh))\bar\lambda(c_jh)}\\
&=& \nabla g(u(c_ih))^\top M^{-1}\left[ p_0 + h\sum_{j=0}^{s-1}\int_0^{c_i}P_j(x)\dd x\,\psi_j(u) + O(h^{s+1})\right], ~ i=1,\dots,s.\nonumber
\end{eqnarray}
Subtracting, side by side,  (\ref{azzera1}) from (\ref{azzera2}) then gives:
\begin{eqnarray*} 
\nabla g(u(c_ih))^\top M^{-1}  \sum_{\ell=0}^{s-1} \int_0^{c_i}P_\ell(x)\dd x \sum_{j=1}^s b_j P_\ell(c_j)\nabla g(u(c_jh))\delta\lambda_j
&=& O(h^s),\\   i=1,\dots,s,
\end{eqnarray*}
from which the statement follows. In fact, the previous equations can be rewritten as (see (\ref{lam}))
$$
\left( \nabla \bfg^\top \, \I_s\P_s^\top\Omega\otimes M^{-1}\, \nabla\bfg\right)\delta\bflam = O(h^s), 
$$
having set $\delta\bflam = \pmatrix{ccc}\delta\lambda_1^\top & \dots & \delta\lambda_1^\top\endpmatrix^\top$.\,\QED

\section{Accuracy}\label{order}
Let us now study the accuracy of the approximations (\ref{q1p1}) and (\ref{q1p1def}). For this purpose, let us denote by\,\footnote{We shall use either one of such notations, depending on the needs.}
\begin{equation}\label{sol}
y(t):= \pmatrix{c} q(t)\\ p(t)\endpmatrix \equiv y(t,\xi,\eta),\qquad \eta =\pmatrix{c}\eta_1\\ \eta_2\endpmatrix\in\RR^{2m},
\end{equation}
the solution of the problem
\begin{equation}\label{local}
\dot y \equiv \pmatrix{c} \dot q \\ \dot p \endpmatrix =
\pmatrix{c} M^{-1}p\\  \nabla U(q)-\nabla g(q)\lambda\endpmatrix =: f(y), \qquad t>\xi, \qquad y(\xi) = \eta,
\end{equation}
where we assume that, according to (\ref{lamq}) $\lambda(t)=\lambda(q(t),p(t))$ is a known function, given by (\ref{inteq}). 
Accordingly, we have denoted by $\bar\lambda(t)$ that defined in (\ref{blam}). We here report standard perturbation results for the solution of (\ref{local}), w.r.t. its parameters (see, e.g., \cite{HNW2008,LIMbook2016}).

\begin{lem}\label{pert}
With reference to (\ref{sol})-(\ref{local}), one has:
$$\frac{\partial}{\partial t} y(t) = f(y(t)), \qquad \frac{\partial}{\partial \xi} y(t) = -\Phi(t,\xi,\eta)f(\eta),$$
where $$\Phi(t,\xi,\eta) = \frac{\partial}{\partial \eta} y(t) \in\RR^{2m\times 2m}.$$
\end{lem}

\medskip
We also need the expansion, with reference to (\ref{blam}),
\begin{eqnarray}\label{Gamj}
\nabla g(u(ch))\bar\lambda(ch) &=& \sum_{j\ge0} P_j(c) \Gamma_j, \qquad c\in[0,1],\\ \nonumber
\Gamma_j &=& \int_0^1 P_j(\tau)\nabla g(u(\tau h))\bar\lambda(\tau h)\dd\tau, \qquad j=0,1,\dots,
\end{eqnarray}
and the following result.

\begin{lem}\label{simme}
The approximation procedure (\ref{gampsij})--(\ref{q1p1def}) and (\ref{azzera0}) is symmetric. Therefore, its convergence order is even.
\end{lem}
\proof
The first part of the statement can be proved by using arguments similar to those used in \cite[Section\,3.3.4]{LIMbook2016}, taking into account that the Gauss-Legendre abscissae $\{c_1,\dots,c_s\}$ are symmetrically distributed in the interval $[0,1]$. The second part then follows from \cite[Theorem\,3.2]{HLW2006}.\,\QED\bigskip

We can now discuss the accuracy of the method, following the framework defined in \cite{BIT2012}.

\begin{theo}\label{ord2s}
With reference to (\ref{gampsij})--(\ref{q1p1def}) and (\ref{azzera0}), one has that the approximation procedure has convergence order $r$, i.e., 
$$q(h)-q_1=O(h^{r+1}), \qquad p(h)-p_1=O(h^{r+1}),$$ 
where $r=s$ if $s$ is even, or $r=s+1$, when $s$ is odd.  In case (see (\ref{blam}) and (\ref{dli})) 
\begin{equation}\label{dli0}
\delta\lambda_i=0, \qquad i=1,\dots,s, 
\end{equation}
then $r=2s$.
\end{theo}
\proof By using the notation (\ref{sol})-(\ref{local}), the results of Theorem~\ref{dlamj}, Lemmas~\ref{Ghj}--\ref{simme},  (\ref{Gamj}), (\ref{quaderr}), (\ref{lim}), and setting
\begin{equation}\label{appr}
y_i = \pmatrix{c} q_i\\ p_i\endpmatrix,\quad i=0,1,\qquad \sigma(t) = \pmatrix{c} u(t)\\ v(t)\endpmatrix,\quad t\in[0,h],  
\end{equation}
we obtain:
\begin{eqnarray*}\lefteqn{
y(h)-y_1 = y(h,0,y_0)-y(h,h,y_1) = y(h,0,\sigma(0))-y(h,h,\sigma(h))}\\
&=&\int_h^0 \frac{\dd}{\dd t} y(h,t,\sigma(t))\dd t\\
&=& \int_h^0 \left[\left.\frac{\partial}{\partial \xi}y(h,\xi,\sigma(t))\right|_{\xi=t}+\left.\frac{\partial}{\partial \eta}y(h,t,\eta)\right|_{\eta=\sigma(t)}\dot\sigma(t)\right]\dd t\\
&=&\int_0^h\Phi(h,t,\sigma(t))\left[f(\sigma(t))-\dot\sigma(t)\right]\dd t \\
&=& h\int_0^1\Phi(h,ch,\sigma(ch))\left[f(\sigma(ch))-\dot\sigma(ch)\right]\dd c \\
&=& h\int_0^1\Phi(h,ch,\sigma(ch))\left[\sum_{j\ge0} P_j(c) \pmatrix{c} M^{-1}\gamma_j(v)\\ \psi_j(u)-\Gamma_j\endpmatrix \right.\\
&&\left.-\sum_{j=0}^{s-1} P_j(c) \pmatrix{c} M^{-1}\gamma_j(v)\\ \psi_j(u)\endpmatrix  +\sum_{j=1}^s \ell_j(c)\pmatrix{c} 0\\ \rho_j\endpmatrix\right]\dd c ~=:~(*).
\end{eqnarray*}
Considering that, from the result of Theorem~\ref{dlamj},
\begin{eqnarray*}
\lefteqn{\sum_{j=1}^s \ell_j(c)\pmatrix{c} 0\\ \rho_j\endpmatrix \,=\, \sum_{j=1}^s \ell_j(c)\pmatrix{c} 0\\ \nabla g(u(c_jh))\lambda_j\endpmatrix}\\
&=&\sum_{j=0}^{s-1} P_j(c)\pmatrix{c} 0\\ \sum_{i=1}^s b_iP_j(c_i)\nabla g(u(c_ih))\lambda_i \endpmatrix
\,=\,  \sum_{j=0}^{s-1} P_j(c)\pmatrix{c}0\\ \Gamma_j +O(h^s)\endpmatrix,
\end{eqnarray*}
it follows  that
\begin{eqnarray*}
(*) &=&h\underbrace{\sum_{j=0}^{s-1} \overbrace{\int_0^1P_j(c)\Phi(h,ch,\sigma(ch))\dd c}^{=\,O(h^j)}\pmatrix{c} 0\\ O(h^s)\endpmatrix}_{=\,O(h^s)}~+~\\
&&h\overbrace{\sum_{j\ge s} \underbrace{\int_0^1P_j(c)\Phi(h,ch,\sigma(ch))\dd c}_{=\,O(h^j)}\underbrace{\pmatrix{c} M^{-1}\gamma_j(v)\\ \psi_j(u)-\Gamma_j\endpmatrix}_{=\,O(h^j)}}^{=\,O(h^{2s})} ~=~O(h^{s+1}).
\end{eqnarray*}
Consequently, the statement follows by considering that $r=s$ at least, and, indeed, $r=s+1$ when $s$ is odd, according to Lemma~\ref{simme}.

In the case where (\ref{dli0}) holds true, then (see (\ref{Gamj})-(\ref{quaderr}))
$$\sum_{i=1}^s b_iP_j(c_i)\nabla g(u(c_ih))\lambda_i  = \Gamma_j + \Delta_j(h), \qquad j=0,\dots,s-1,$$
so that
$$
\sum_{j=1}^s \ell_j(c)\pmatrix{c} 0\\ \rho_j\endpmatrix =  \sum_{j=0}^{s-1} P_j(c)\pmatrix{c}0\\ \Gamma_j +\Delta_j(h)\endpmatrix,
$$
and, consequently,
\begin{eqnarray*}
(*) &=&h\underbrace{\sum_{j=0}^{s-1} \overbrace{\int_0^1P_j(c)\Phi(h,ch,\sigma(ch))\dd c}^{=\,O(h^j)}\overbrace{\pmatrix{c} 0\\ \Delta_j(h)\endpmatrix}^{=\,O(h^{2s-j})}}_{=\,O(h^{2s})}~+~\\
&&h\overbrace{\sum_{j\ge s} \underbrace{\int_0^1P_j(c)\Phi(h,ch,\sigma(ch))\dd c}_{=\,O(h^j)}\underbrace{\pmatrix{c} M^{-1}\gamma_j(v)\\ \psi_j(u)-\Gamma_j\endpmatrix}_{=\,O(h^j)}}^{=\,O(h^{2s})} ~=~O(h^{2s+1}).\,\QED
\end{eqnarray*}

\begin{cor}\label{norm}
In the hypothesis of the previous Theorem~\ref{ord2s}, and using the notation (\ref{sol}) and (\ref{appr}), one has:\,\footnote{Hereafter, $|\cdot|$ denotes any convenient vector norm.}
$$\|y-\sigma\| := \max_{\tau\in[0,1]} | y(\tau h)-\sigma(\tau h)| = O(h^{s+1}).$$
\end{cor}
\proof In fact, one has:
$$
y(\tau h)-\sigma(\tau h) = y(\tau h,0,\sigma(0))-y(\tau h,\tau h,\sigma(\tau h)) =
\int_{\tau h}^0 \frac{\dd}{\dd t} y(\tau h,t,\sigma(t))\dd t =:(**).
$$
By repeating similar steps as those in the proof of Theorem~\ref{ord2s}, one then obtains:
\begin{eqnarray*}
(**) &=&h\underbrace{\sum_{j=0}^{s-1} \overbrace{\int_0^\tau P_j(c)\Phi(\tau h,ch,\sigma(ch))\dd c}^{=\,O(h^0)}\pmatrix{c} 0\\ O(h^s)\endpmatrix}_{=\,O(h^s)}~+~\\
&&h\overbrace{\sum_{j\ge s} \underbrace{\int_0^\tau P_j(c)\Phi(\tau h,ch,\sigma(ch))\dd c}_{=\,O(h^0)}\underbrace{\pmatrix{c} M^{-1}\gamma_j(v)\\ \psi_j(u)-\Gamma_j\endpmatrix}_{=\,O(h^j)}}^{=\,O(h^s)} ~=~O(h^{s+1}).\,\QED
\end{eqnarray*}

\subsection{Hidden constraints and Lagrange multipliers}\label{acclam}

We can now confirm that the value of $\bar r$ in (\ref{hidcons}) has to be at least equal to the order of the method, as stated in Theorem~\ref{ord2s}.
In fact, at the first step one obtains:
\begin{eqnarray*}
0 &=& \nabla g(q(h))^\top M^{-1} p(h)\\
&=& \left[\nabla g(u(h)) + O(h^{r+1})\right]^\top M^{-1}\left[ v(h) + O(h^{r+1})\right]\\
&=& \nabla g(u(h))^\top M^{-1} v(h) + O(h^{r+1}).
\end{eqnarray*}
As is clear, at the generic $n$th step, setting $t_n=nh$, and
\begin{eqnarray}\nonumber
u_n(ch) &\approx& q_n(ch):=q(t_n+ch), \\ \label{unvn}
v_n(ch)&\approx& p_n(ch) :=p(t_n+ch), \qquad c\in[0,1],
\end{eqnarray}
the polynomial approximations and the restriction of the solution in the interval $[t_n,t_{n+1}]$, 
since $u_n(h) = q_n(h)+O(h^r)$ and $v_n(h)=p_n(h)+O(h^r)$, one has:
$$ \nabla g(u_n(h))^\top M^{-1} v_n(h) = O(h^r).$$

Similarly, by using the same notation above, we can discuss the accuracy of the approximate Lagrange multipliers, 
\begin{equation}\label{ln}
\lambda_i^n \approx \lambda_n(c_ih) := \lambda(t_n+c_ih), \qquad i=1,\dots,s,
\end{equation}
The following result holds true.

\begin{theo}\label{ordlam}
With reference to (\ref{gampsij})--(\ref{q1p1def}) and (\ref{azzera0}), assume that the procedure 
has convergence order $r$, as stated in Theorem~\ref{ord2s}, in the case where (\ref{dli0}) does not hold.
Then, the Langrange multipliers are approximated with at least order $r-1$.
\end{theo}
\proof We recall that, in the considered case, the order of the method is $r=s$, if $s$ is even, or $r=s+1$, otherwise. With reference to 
(\ref{unvn})-(\ref{ln}), the approximate Lagrange multipliers then satisfy:
\begin{eqnarray}\label{errln} \lefteqn{
h \nabla g(u_n(c_ih))^\top M^{-1}  \sum_{\ell=0}^{s-1} \int_0^{c_i}P_\ell(x)\dd x \sum_{j=1}^s b_j P_\ell(c_j)\nabla g(u_n(c_jh))\lambda_j^n}\\  
&=& \nabla g(u_n(c_ih))^\top M^{-1}\left[ v_n(0) + h\sum_{j=0}^{s-1}\int_0^{c_i}P_j(x)\dd x\,\psi_j(u_n)\right], \quad i=1,\dots,s.\nonumber
\end{eqnarray}
On the other hand, from (\ref{inteq}) one deduces that the exact multipliers satisfy:
\begin{eqnarray*}\nonumber\lefteqn{
h \nabla g(q_n(c_ih))^\top M^{-1}  \sum_{\ell\ge 0} \int_0^{c_i}P_\ell(x)\dd x \int_0^1 P_\ell(c) \nabla g(q_n(ch))\lambda_n(ch)\dd c}\\  
&=& \nabla g(q_n(c_ih))^\top M^{-1}\left[ p_n(0) + h\sum_{j\ge 0}\int_0^{c_i}P_j(x)\dd x\,\psi_j(q_n)\right], \quad i=1,\dots,s.\nonumber
\end{eqnarray*}
Considering that, by virtue of Corollary~\ref{norm},
$$\|u_n-q_n\| = \|v_n-p_n\| = O(h^r) + O(h^{s+1}) = O(h^r),$$
from the last equation one derives:
\begin{eqnarray*}\nonumber\lefteqn{
h \nabla g(u_n(c_ih))^\top M^{-1}  \sum_{\ell= 0}^{s-1} \int_0^{c_i}P_\ell(x)\dd x \int_0^1 P_\ell(ch) \nabla g(u_n(ch))\lambda_n(ch)\dd c}\\  
&=& \nabla g(u_n(c_ih))^\top M^{-1}\left[ p_n(0) + h\sum_{j=0}^{s-1}\int_0^{c_i}P_j(x)\dd x\,\psi_j(u_n)\right] + O(h^r), \\
&& \qquad i=1,\dots,s,\nonumber
\end{eqnarray*}
and, therefore,
\begin{eqnarray}\nonumber \lefteqn{
h \nabla g(u_n(c_ih))^\top M^{-1}  \sum_{\ell=0}^{s-1} \int_0^{c_i}P_\ell(x)\dd x \sum_{j=1}^s b_j P_\ell(c_j)\nabla g(u_n(c_jh))\lambda_n(c_jh)}\\  \nonumber
&=& \nabla g(u_n(c_ih))^\top M^{-1}\left[ v_n(0) + h\sum_{j=0}^{s-1}\int_0^{c_i}P_j(x)\dd x\,\psi_j(u_n)\right] + O(h^r), \\
&&\qquad i=1,\dots,s. \label{errln1}
\end{eqnarray}
From (\ref{errln}) and (\ref{errln1}), one eventually obtains:
\begin{eqnarray*}\nonumber \lefteqn{
\nabla g(u_n(c_ih))^\top M^{-1}  \sum_{\ell=0}^{s-1} \int_0^{c_i}P_\ell(x)\dd x \sum_{j=1}^s b_j P_\ell(c_j)\nabla g(u_n(c_jh))(\lambda_n(c_jh)-\lambda_j^n)}\\  \nonumber
&=& O(h^{r-1}), \qquad i=1,\dots,s,\hspace{7cm}
\end{eqnarray*}
from which the statement follows.\,\QED\bigskip

\begin{rem}
Clearly, in the more favourable case where (\ref{dli0}) holds true, the accuracy of the approximate Lagrange multipliers could be better.
\end{rem} 

\section{Discretization}\label{discrete}

For the actual implementation of the approximation procedure (\ref{gampsij})--(\ref{q1p1def}) let us recast it into a more operative form:
\begin{eqnarray}\label{lim1}
u(c_ih) &=& q_0 + hM^{-1}\sum_{j=0}^{s-1} \int_0^{c_i} P_j(x)\dd x\,\gamma_j(v), \\ \nonumber
v(c_ih) &=& p_0 + h\sum_{j=0}^{s-1} \int_0^{c_i} P_j(x)\dd x\,\left[ \psi_j(u)- \zeta_j(u,\bflam)\right],\qquad i=1,\dots,s,
\end{eqnarray}
where, for $j=0,\dots,s-1$,
\begin{eqnarray}\nonumber
\gamma_j(v) &=& \int_0^1 P_j(c)v(ch)\dd c ~\equiv~ \sum_{i=1}^s b_i P_j(c_i)v(c_ih),\\ \label{lim2}
\psi_j(u)        &=& \int_0^1 P_j(c)\nabla U(u(ch))\dd c,\\ \nonumber
\zeta_j(u,\bflam) &=& \int_0^1 P_j(c)\sum_{i=1}^s\ell_j(c)\nabla g(u(c_ih))\lambda_i ~\equiv~ \sum_{i=1}^s b_iP_j(c_i)\nabla g(u(c_ih))\lambda_i,
\end{eqnarray}
and the vector $\bflam$, defined in (\ref{vettori}), is the solution of (\ref{lam}). Moreover, should $U$ be a quadratic potential, then one could exactly compute
$$\psi_j(u) = \sum_{i=1}^s b_i P_j(c_i)\nabla U(u(c_ih)),\qquad j=0,\dots,s-1,$$
but, in general, this could be not the case: as matter of fact, many potentials are rational functions. Consequently, we need to use a more accurate quadrature, for approximating the Fourier coefficients $\psi_j(u)$, $j=0,\dots,s-1$. For this purpose, following the same approach used for HBVM$(k,s)$ methods (see, e.g., \cite{LIMbook2016,axioms2018}) we can use a Gauss-Legendre quadrature of order $2k\ge 2s$, with weights and abscissae
$$(\hat c_i, \hat b_i), \qquad i=1,\dots,k,$$
thus obtaining:
\begin{eqnarray}\label{psij}
\psi_j(u) &=& \sum_{i=1}^k \hat b_i P_j(\hat c_i)\nabla U(u(\hat c_ih)) \,+\,\hat\Delta_j(h)\\
              &=:& \hat\psi_j(u) \,+\,\hat\Delta_j(h) ,\qquad j=0,\dots,s-1,\nonumber
\end{eqnarray}
with
\begin{equation}\label{hDj}
\hat\Delta_j(h) = \left\{\begin{array}{cc} 0, & ~if~ U\in\Pi_\mu, ~with~ \mu\le 2k/s, \\[2mm]
O(h^{2k-j}), &otherwise.\end{array}\right.
\end{equation}

\begin{rem}\label{costo}
We shall see that we can choose $k$ as large as it will be needed, without increasing the computational cost of the method significantly.
Consequently, the quadrature errors can be made arbitrarily small, until they fall within the round-off error level of the used finite precision arithmetic.
\end{rem}

The following result clearly holds true.
\begin{theo}\label{ovvio}
In case $U\in\Pi_\mu$, with $\mu\le 2k/s$, then $\psi_j(u)=\hat\psi_j(u)$, $j=0,\dots,s-1$,  and the results of Theorems~\ref{gheq0}, \ref{evai}, \ref{ord2s}, and \ref{ordlam}, as well as Corollary~\ref{norm}, continue to hold.
\end{theo}
 
Conversely, by continuing to denote by $u$ and $v$ the new polynomial approximations using $\hat\psi_j$ defined by (\ref{psij}), instead of $\psi_j$ defined in (\ref{gampsij}), the statements of the previous results modify as follows, the proofs being similar.

\begin{theo}\label{evaik} Assume that:  
\begin{equation}\label{azzera}
\nabla g(u(c_ih))^\top M^{-1} v(c_ih)=0\in\RR^\nu, \qquad i=1,\dots,s.
\end{equation}
Then, for all $k\ge s$, properties (\ref{Hgcons}) are replaced by 
\begin{eqnarray*}
H(q_1,p_1) &=& \left\{ \begin{array}{cc}  H(q_0,p_0), &\mbox{if~} U\in\Pi_\mu, \mbox{~with~} \mu\le 2k/s,\\[2mm]
H(q_0,p_0) + O(h^{2k+1}), &\mbox{otherwise},\end{array}\right. \\[2mm]
g(q_1)&=&0,
\end{eqnarray*}
whereas for (\ref{hidcons}) the considerations made in Section~\ref{acclam} continue to hold.
\end{theo}

\begin{theo}\label{ord2sk}
With reference to (\ref{gampsij})--(\ref{q1p1def}) and (\ref{azzera}), with $\psi_j(u)$ formally replaced by $\hat\psi_j(u)$ as defined in (\ref{psij}), and assuming $k\ge s$, the approximation procedure has convergence order $r$, i.e., 
$$q(h)-q_1=O(h^{r+1}), \qquad p(h)-p_1=O(h^{r+1}),$$ 
where $r=s$ if $s$ is even, or $r=s+1$, when $s$ is odd. In the case where (\ref{dli0}) holds true, then $r=2s$.
\end{theo}

\begin{cor}\label{normk}
In the hypothesis of the previous Theorem~\ref{ord2sk}, and using the notation (\ref{sol}) and (\ref{appr}), one has:
$$\|y-\sigma\| := \max_{\tau\in[0,1]} | y(\tau h)-\sigma(\tau h)| = O(h^{s+1}).$$
\end{cor}

\begin{theo}\label{ordlamk}
In the hypothesis of the previous Theorem~\ref{ord2sk}, the result of Theorem~\ref{ordlam} continues to hold.\end{theo}

We now prove Theorem~\ref{evaik}, since the other statements can be similarly proved, by slightly modifying the proof of Theorems~\ref{ord2s} and \ref{ordlam}, and Corollary~\ref{norm}, respectively (the procedure is similar to that used in \cite{BIT2012} when discretizing the involved integrals).

\bigskip
\proof (of Theorem~\ref{evaik}). We prove only the statement concerning the Hamiltonian, since the remaining ones are straightforward. 
In the polynomial case, the statement follows from Theorem~\ref{ovvio}. Conversely, for all $k\ge s$, by taking into account (\ref{psij})-(\ref{hDj}), one has:
\begin{eqnarray*}\nonumber
\lefteqn{H(q_1,p_1)-H(q_0,p_0) ~\equiv~ H(u(h),v(h))-H(u(0),v(0))}\\ \nonumber
&=& \int_0^h \frac{\dd}{\dd t} H(u(t),v(t))\dd t\\ \nonumber
&=& h\int_0^1 \left[H_q(u(ch),v(ch))^\top \dot u(ch)+ H_p(u(ch),v(ch))^\top \dot v(ch)\right]\dd c\\ \nonumber
&=& h\int_0^1\left[-\nabla U(u(ch))^\top M^{-1}\sum_{j=0}^{s-1}P_j(c)\gamma_j(v)+v(ch)^\top M^{-1}\sum_{j=0}^{s-1}P_j(c)\hat\psi_j(u)\right.\\ \nonumber
&&\left.\qquad ~~-~v(ch)^\top M^{-1}\sum_{j=1}^s\ell_j(c)\rho_j\right]\dd c\\ \nonumber
&=& -h \sum_{j=0}^{s-1} \underbrace{\left[\int_0^1P_j(c)\nabla U(u(ch))\dd c\right]^\top}_{=\,\psi_j(u)^\top}M^{-1}\gamma_j(v)\\ \nonumber
&& +h \sum_{j=0}^{s-1} \underbrace{\left[\int_0^1P_j(c)v(ch)\dd c\right]^\top}_{=\,\gamma_j(u)^\top} M^{-1} \hat\psi_j(u) 
-h\int_0^1 \sum_{j=1}^s v(ch)^\top M^{-1} \ell_j(c)\rho_j\dd c\\ \nonumber
&=& -h\sum_{j=0}^{s-1} \psi_j(u)^\top M^{-1} \gamma_j(v)+h\sum_{j=0}^{s-1}\gamma_j(v)^\top M^{-1}\underbrace{[\psi_j(u)~-\overbrace{\hat\Delta_j(h)}^{=\,O(h^{2k-j})}]}_{=\,\hat\psi_j(u)}\\
&&-h\sum_{j=1}^s b_j v(c_jh)^\top M^{-1}\rho_j\\ 
&=& -h\sum_{j=0}^{s-1}\underbrace{\gamma_j(v)^\top M^{-1} \hat\Delta_j(h)}_{=\,O(h^{2k})}\,-h\sum_{j=1}^s b_j \underbrace{v(c_jh)^\top M^{-1} \nabla g(u(c_jh))}_{=\,0}\lambda_j \\
&=& O(h^{2k+1}).\,\QED
\end{eqnarray*}

\subsection{Discrete problem in vector form}\label{dispro}
Let us now derive a vector formulation of the discrete problem obtained by formally replacing $\psi_j(u)$ in (\ref{lim2}) with the quadrature $\hat\psi_j(u)$ defined in (\ref{psij}). For this purpose, besides the matrices and vectors defined in (\ref{matrici}) and (\ref{vettori}), we also need the following ones:
\begin{eqnarray}\nonumber
\hat\P_s &=&\pmatrix{ccc} P_0(\hat c_1)&\dots&P_{s-1}(\hat c_1)\\ \vdots & &\vdots\\ P_0(\hat c_k)& \dots &P_{s-1}(\hat c_k)\endpmatrix\in\RR^{k\times s},\\[2mm] \nonumber
\hat\I_s &=&\pmatrix{ccc} \int_0^{\hat c_1}P_0(x)\dd x&\dots&\int_0^{\hat c_1}P_{s-1}(x)\dd x\\ \vdots & &\vdots\\ \int_0^{\hat c_k}P_0(x)\dd x& \dots &\int_0^{\hat c_k}P_{s-1}(x)\dd x\endpmatrix\in\RR^{k\times s},\\[2mm]  \label{new}
\hat\Omega &=&\diag(\hat b_1,\,\dots,\,\hat b_k)\in\RR^{k\times k},\\[2mm] \nonumber
\hat{\bfpsi} &=& \pmatrix{c} \hat\psi_0(u)\\ \vdots\\ \hat\psi_{s-1}(u)\endpmatrix, ~
\bfgamma \,=\, \pmatrix{c} \gamma_0(v)\\ \vdots\\ \gamma_{s-1}(v)\endpmatrix \in\RR^{sm},\\[2mm] \nonumber
\bfzeta &=& \pmatrix{c} \zeta_0(u,\bflam)\\ \vdots\\ \zeta_{s-1}(u,\bflam)\endpmatrix \in\RR^{sm},\qquad 
\hat\bfe     \,=\, \pmatrix{c}1\\ \vdots\\ 1\endpmatrix\in\RR^k.
\end{eqnarray}

One then obtains, by also taking into account (\ref{lam}), the following discrete problem: 
\begin{eqnarray}\nonumber
\nabla\bfg &=& \nabla g(\bfe\otimes q_0 + h\I_s\otimes M^{-1} \bfgamma)\\ \nonumber
\bflam       &=& \left( \nabla \bfg^\top \, h\I_s\P_s^\top\Omega\otimes M^{-1}\, \nabla\bfg\right)^{-1}
\nabla\bfg^\top \left[  \bfe\otimes M^{-1}p_0 + h\I_s\otimes M^{-1}\hat\bfpsi\right]\\  \label{dispro1}
\bfzeta      &=& \P_s^\top\Omega\otimes I_m (\nabla\bfg\cdot\bflam)  \\ \nonumber
\bfgamma &=& \P_s^\top\Omega\otimes I_m \left[ \bfe\otimes p_0 +h\I_s\otimes I_m (\hat\bfpsi -\bfzeta)\right]\\ \nonumber 
\hat\bfpsi  &=&  \hat\P_s^\top \hat\Omega\otimes I_m \, \nabla U(\hat\bfe\otimes q_0 + h\hat\I_s\otimes M^{-1} \bfgamma) 
\end{eqnarray}
where, considering that 
\begin{eqnarray*}
\bfe\otimes q_0 + h\I_s\otimes M^{-1} \bfgamma &\equiv& \pmatrix{c} u(c_1 h)\\ \vdots \\ u(c_s h)\endpmatrix,\\
\hat\bfe\otimes q_0 + h\hat\I_s\otimes M^{-1} \bfgamma &\equiv& \pmatrix{c} u(\hat c_1 h)\\ \vdots \\ u(\hat c_k h)\endpmatrix,\\
\bfe\otimes p_0 +h\I_s\otimes I_m(\hat\bfpsi -\bfzeta) &\equiv&  \pmatrix{c} v(c_1 h)\\ \vdots \\ v(c_s h)\endpmatrix,
\end{eqnarray*}
it is meant that:
$$\nabla U(\hat\bfe\otimes q_0 + h\hat\I_s\otimes M^{-1} \bfgamma) = \pmatrix{c} \nabla U(u(\hat c_1 h))\\ \vdots \\ \nabla U(u(\hat c_k h))\endpmatrix\in\RR^{km},$$
and, similarly,
$\nabla g(\bfe\otimes q_0 + h\I_s\otimes M^{-1} \bfgamma)$ computes the $sm\times s\nu$ matrix $\nabla\bfg$ defined in (\ref{matrici}).  

\begin{rem} In the case where the holonomic constraints are not present, i.e. $\bfzeta=\bfzero$ in (\ref{dispro1}), the first three equations are void. Moreover,
the last two equations, combined together, provide the discrete problem of a HBVM$(k,s)$ method applied to the corresponding unconstrained problem.
For this reason, we shall continue to refer to (\ref{dispro1}) as a HBVM$(k,s)$ method applied for solving the constrained Hamiltonian problem (\ref{constrp}).
\end{rem}

We conclude this section by observing that the discrete problem (\ref{dispro1}), which has to be solved within machine accuracy, in order to gain the conservation properties, can be numerically solved by suitably modifying the Newton-type iterations studied in \cite{BIT2011,BFCI2014} (see also \cite{BrMa2002,BrMa2009}). However, for carrying out the numerical tests  of Section~\ref{numex}, we shall use the following straightforward fixed-point iteration,
\begin{eqnarray}\nonumber
\nabla\bfg^\ell &=& \nabla g(\bfe\otimes q_0 + h\I_s\otimes M^{-1} \bfgamma^\ell)\\  \nonumber
\bflam^{\ell} &=& \left( (\nabla \bfg^\ell)^\top \, h\I_s\P_s^\top\Omega\otimes M^{-1}\, \nabla\bfg^\ell\right)^{-1}\\ \nonumber
&& \qquad (\nabla\bfg^\ell)^\top \left[  \bfe\otimes M^{-1}p_0 + h\I_s\otimes M^{-1}\hat\bfpsi^\ell\right] \\ \label{iter}
\bfzeta^{\ell}      &=& \P_s^\top\Omega\otimes I_m (\nabla\bfg^\ell\cdot\bflam^{\ell}) \\ \nonumber
\bfgamma^{\ell+1} &=& \P_s^\top\Omega\otimes I_m \left[ \bfe\otimes p_0 +h\I_s\otimes I_m (\hat\bfpsi^\ell -\bfzeta^\ell)\right]\\ \nonumber
\hat\bfpsi^{\ell+1}  &=&  \hat\P_s^\top \hat\Omega\otimes I_m \, \nabla U(\hat\bfe\otimes q_0 + h\hat\I_s\otimes M^{-1} \bfgamma^\ell),\qquad \quad\ell=0,1,\dots,
\end{eqnarray}
which can be conveniently started from
$$\bfgamma^0=\hat\bfpsi^0=\bfzero\in\RR^{sm}.$$
Assuming a continuous Lipschitz condition for $\nabla U$, and recalling that the constraints are quadratic, the convergence of the previous iteration can be seen to be granted for all sufficiently small timesteps $h$. 
Once (\ref{iter}) has been iterated until convergence, the new approximations are given by:
\begin{equation}\label{q1p1new}
q_1 = q_0 + h\gamma_0, \qquad p_1 = p_0+h(\hat\psi_0 -\zeta_0), 
\end{equation}
where we have removed the arguments of the Fourier coefficients, in order to have a more algorithmic description. An approximation of the vector of the Lagrange multipliers can be similarly obtained as
\begin{equation}\label{lagmul}
\bar\lambda = \sum_{i=1}^s \ell_i(1)\lambda_i,
\end{equation}
$\lambda_1,\dots,\lambda_s$ being, according to (\ref{vettori}), the (block) entries of the vector $\bflam$ at convergence.

\section{Numerical examples}\label{numex}
In this section we report a number of test problems, which we solve by means of HBVM$(k,s)$ methods. The methods have been implemented in Matlab (Rel.~2020b) on a 3GHz Intel Xeon W10 core computer with 64GB RAM. In all numerical tests, the discrete problem (\ref{dispro1}) is
solved according to (\ref{iter}), which is iterated until full machine accuracy is gained. The different problems considered here, are aimed at assessing different facets of the theoretical results.

\subsection{The simple pendulum}\label{simpend}

The first  example that we consider is the simple pendulum in cartesian coordinates. We place the reference system, so that the motion is in the plane $(q_1,q_2)$. Assuming to have normalized for the mass and the acceleration of gravity, and also assuming a unit length of the pendulum, the resulting equations are (hereafter, $\|\cdot\|$ will denote the Euclidean norm):
\begin{equation}\label{pend}
\dot q = p, \qquad \dot p = -e_2, \qquad g(q) := \|q\|^2-1 = 0, 
\end{equation}
$e_2$ being the second unit vector in $\RR^2$, with Hamiltonian
\begin{equation}\label{pendH}
H(q,p) = \frac{1}2 \|p\|^2 +e_2^\top q.
\end{equation}
Having set the origin at the suspension point of the pendulum, a set of consistent initial conditions are:
\begin{equation}\label{pend0}
q(0) = \pmatrix{c} 0\\ -1\endpmatrix, \qquad p(0) = \pmatrix{c}1\\ 0\endpmatrix.
\end{equation}

\begin{table}[p]
\caption{obtained results for problem (\ref{pend})--(\ref{pend0}) by using a HBVM$(s,s)$ method with timestep (\ref{hi}).}\label{pendtab}

\smallskip\small
\centerline{
\begin{tabular}{|c|cccccccc|}
\hline
\multicolumn{9}{|c|}{$s=1$}\\
\hline
$i$ & $e_y$ & rate & $e_\lambda$ & rate & $e_{hid}$ & rate & $e_g$ & $e_H$\\
\hline
  0 & 5.87e-01 &     --- & 5.47e-01 &    --- & 6.08e-02 &    --- & 5.55e-15 & 2.22e-15  \\ 
  1 & 1.66e-01 &   1.8 & 1.65e-01 &   1.7 & 1.48e-02 &   2.0 & 2.22e-15 & 7.77e-16  \\ 
  2 & 4.62e-02 &   1.8 & 3.58e-02 &   2.2 & 3.67e-03 &   2.0 & 2.44e-15 & 1.11e-15  \\ 
  3 & 1.17e-02 &   2.0 & 1.01e-02 &   1.8 & 9.16e-04 &   2.0 & 1.89e-15 & 1.33e-15  \\ 
  4 & 2.95e-03 &   2.0 & 3.67e-03 &   1.5 & 2.29e-04 &   2.0 & 3.33e-15 & 5.55e-16  \\ 
  5 & 7.37e-04 &   2.0 & 1.56e-03 &   1.2 & 5.72e-05 &   2.0 & 1.33e-15 & 8.88e-16  \\ 
  6 & 1.84e-04 &   2.0 & 7.20e-04 &   1.1 & 1.43e-05 &   2.0 & 3.77e-15 & 2.22e-15  \\ 
  7 & 4.61e-05 &   2.0 & 3.46e-04 &   1.1 & 3.58e-06 &   2.0 & 3.11e-15 & 1.33e-15  \\ 
  8 & 1.15e-05 &   2.0 & 1.69e-04 &   1.0 & 8.94e-07 &   2.0 & 7.11e-15 & 4.66e-15  \\ 
\hline
\hline
\multicolumn{9}{|c|}{$s=2$}\\
\hline
$i$ & $e_y$ & rate & $e_\lambda$ & rate & $e_{hid}$ & rate & $e_g$ & $e_H$\\
\hline
  0 & 4.61e-02 &   --- & 1.14e-02 &   --- & 1.98e-02 &   --- & 8.88e-16 &   5.00e-16  \\ 
  1 & 4.30e-03 &   3.4 & 8.41e-04 &   3.8 & 4.90e-03 &   2.0 & 5.55e-16 &   5.55e-16  \\ 
  2 & 5.99e-04 &   2.8 & 1.75e-03 &  -1.1 & 1.22e-03 &   2.0 & 1.55e-15 &   5.55e-16  \\ 
  3 & 1.20e-04 &   2.3 & 1.31e-03 &   0.4 & 3.05e-04 &   2.0 & 1.33e-15 &   2.78e-16  \\ 
  4 & 2.81e-05 &   2.1 & 7.68e-04 &   0.8 & 7.63e-05 &   2.0 & 8.88e-16 &   4.44e-16  \\ 
  5 & 6.91e-06 &   2.0 & 4.13e-04 &   0.9 & 1.91e-05 &   2.0 & 1.33e-15 &   2.16e-15  \\ 
  6 & 1.72e-06 &   2.0 & 2.14e-04 &   0.9 & 4.77e-06 &   2.0 & 4.33e-15 &   1.22e-15  \\ 
  7 & 4.29e-07 &   2.0 & 1.09e-04 &   1.0 & 1.19e-06 &   2.0 & 3.33e-15 &   1.72e-15  \\ 
  8 & 1.07e-07 &   2.0 & 5.48e-05 &   1.0 & 2.98e-07 &   2.0 & 6.66e-15 &   3.05e-15  \\ 
\hline
\hline
\multicolumn{9}{|c|}{$s=3$}\\
\hline
$i$ & $e_y$ & rate & $e_\lambda$ & rate & $e_{hid}$ & rate & $e_g$ & $e_H$\\
\hline
  0 & 2.83e-03 &   --- & 5.58e-03 &   --- & 1.26e-03 &  --- & 4.44e-16 &   2.22e-16  \\ 
  1 & 1.86e-04 &   3.9 & 2.21e-03 &   1.3 & 6.50e-05 &   4.3 & 4.44e-16 &   2.22e-16  \\ 
  2 & 9.45e-06 &   4.3 & 1.94e-04 &   3.5 & 3.95e-06 &   4.0 & 4.44e-16 &   1.67e-16  \\ 
  3 & 5.86e-07 &   4.0 & 2.40e-05 &   3.0 & 2.46e-07 &   4.0 & 8.88e-16 &   5.55e-16  \\ 
\hline
\hline
\multicolumn{9}{|c|}{$s=4$}\\
\hline
$i$ & $e_y$ & rate & $e_\lambda$ & rate & $e_{hid}$ & rate & $e_g$ & $e_H$\\
\hline
  0 & 3.70e-04 &   --- & 2.89e-03 &   --- & 1.60e-04 &   --- & 2.22e-16 &  1.11e-16  \\ 
  1 & 2.24e-05 &   4.0 & 4.58e-04 &   2.7 & 9.63e-06 &   4.1 & 2.22e-16 &  3.33e-16  \\ 
  2 & 1.44e-06 &   4.0 & 6.04e-05 &   2.9 & 5.93e-07 &   4.0 & 5.55e-16 &   2.22e-16  \\ 
  3 & 9.08e-08 &   4.0 & 7.32e-06 &   3.0 & 3.73e-08 &   4.0 & 1.11e-15 &   5.55e-16  \\ 
\hline
\hline 
\end{tabular}}
\end{table}

Since the Hamiltonian (\ref{pendH}) is quadratic, a HBVM$(s,s)$ method (i.e., the $s$-stage Gauss method) is energy-conserving of order $r=2\lceil\frac{s}2\rceil$, with the 
approximation of the Lagrange multipliers (\ref{lagmul}) of order $r-1$, for all $s=1,2,\dots$. Also the hidden constraints are conserved with order $r$.
In Table~\ref{pendtab} we list the obtained results, by using the HBVM$(s,s)$ method, $s=1,2,3,4$, with timestep 
\begin{equation}\label{hi}
h_i=2^{-i}, \qquad i=0,1,\dots, 
\end{equation}
for solving problem (\ref{pend})--(\ref{pend0}) on the interval $[0,10]$. The error $e_y$ in the solution is numerically estimated. Moreover, we have denoted
by $e_\lambda$ the error in the Lagrange multipliers, $e_{hid}$ the error in the hidden constraints, $e_g$ the error in the constraints, and
by $e_H$ the Hamiltonian error (the same notation will be used in the subsequent tests). As one may see, all the obtained results exactly match the expected ones.

\begin{figure}[t]
\centerline{
\includegraphics[width=6cm]{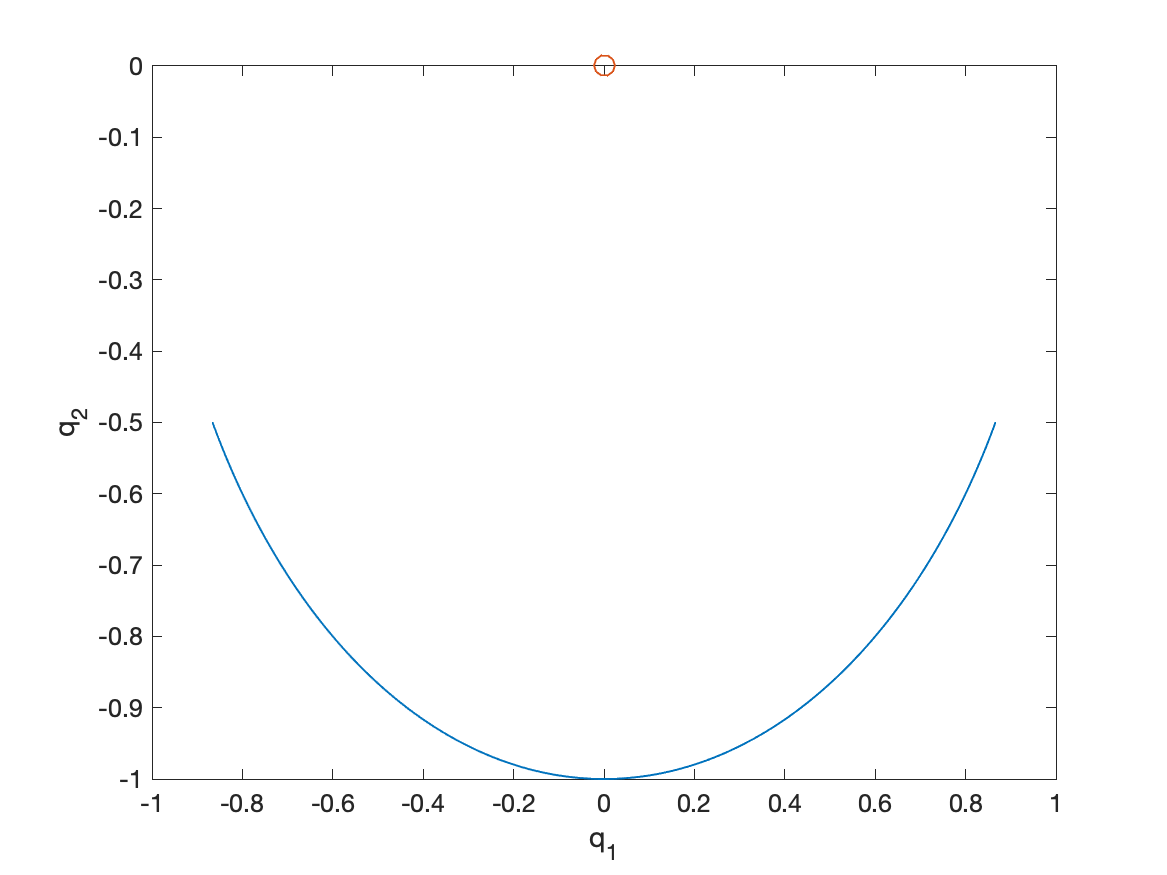}\quad \includegraphics[width=6cm]{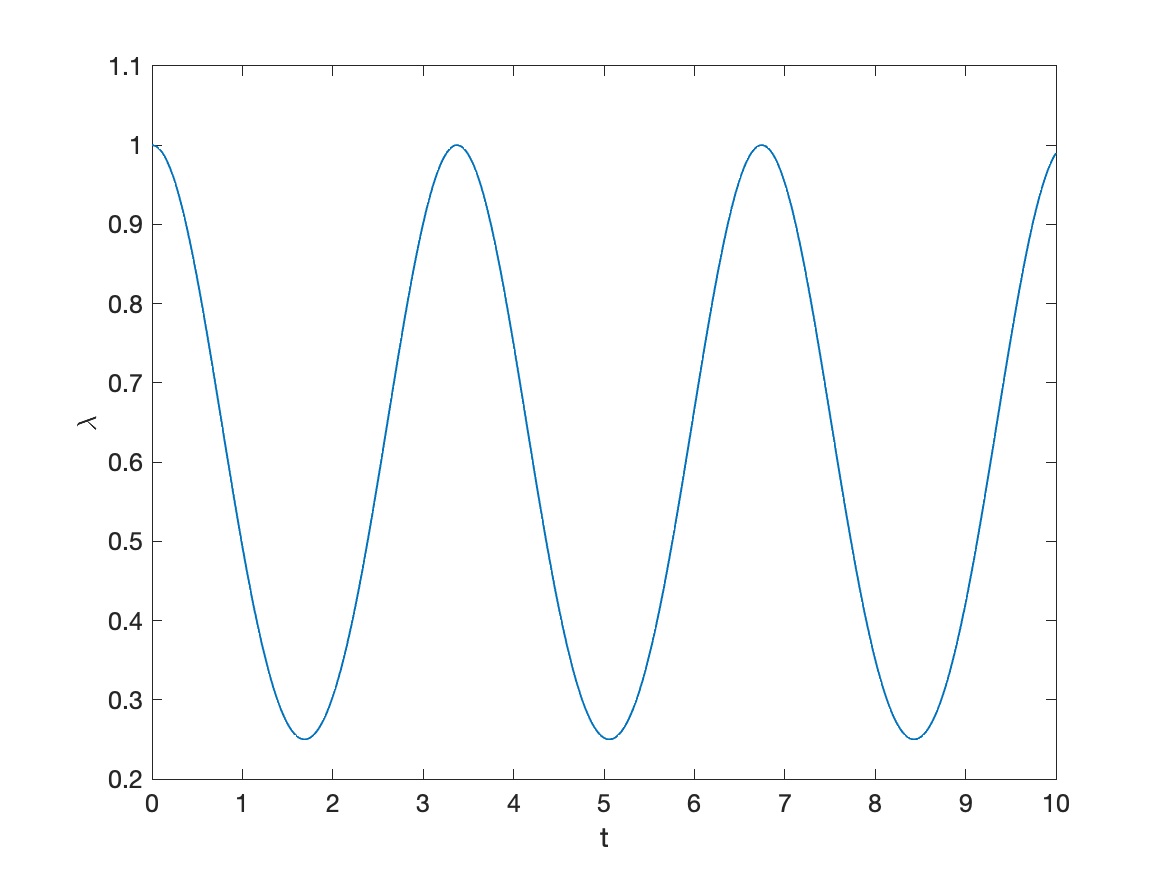}}
\caption{\label{pend_lam} phase portrait (left) and Lagrange multiplier (right) for problem (\ref{pend})--(\ref{pend0}). }
\end{figure}

As is well known, in such a case, the Lagrange multiplier has the physical meaning of the tension of the suspension rod.
Its plot is in Figure~\ref{pend_lam}, together with the phase portrait of the solution.

\subsection{A modified pendulum}\label{modpend} 

\begin{figure}[t]
\centerline{
\includegraphics[width=6cm]{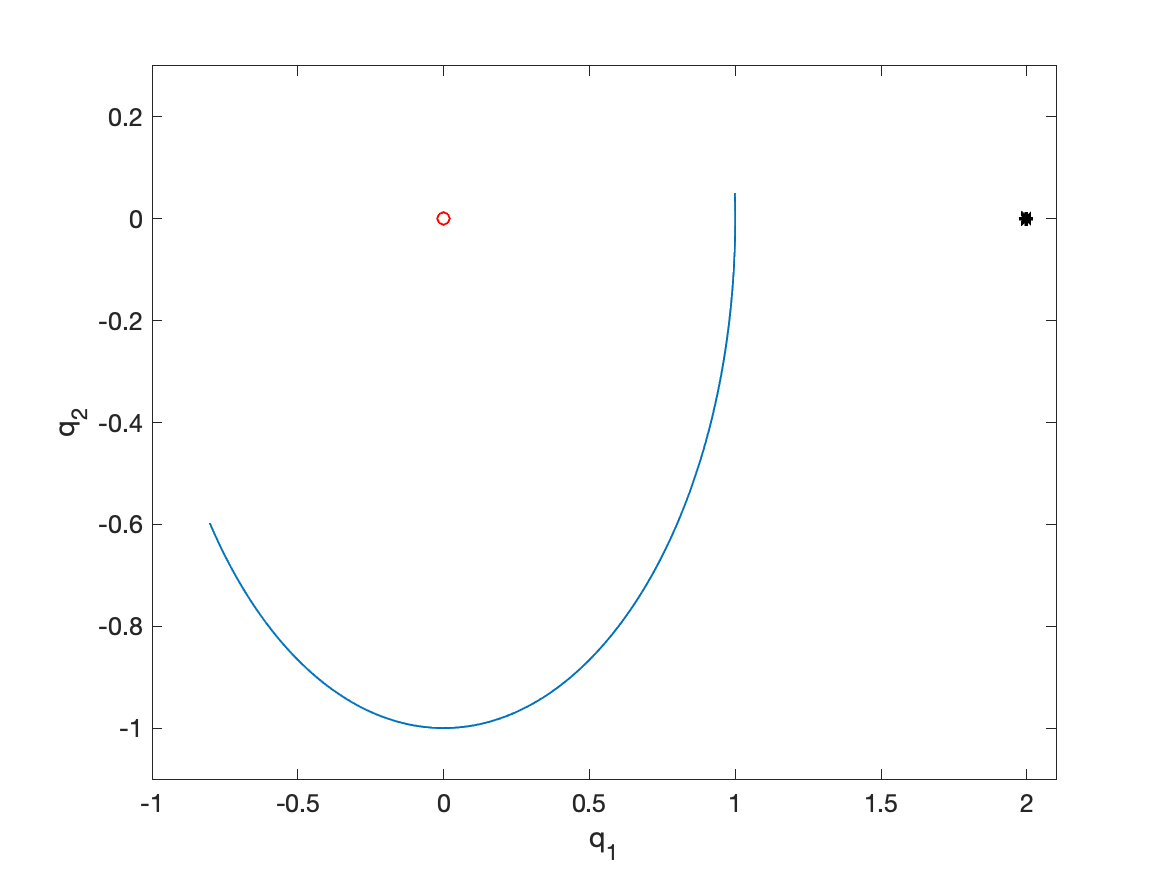}\quad \includegraphics[width=6cm]{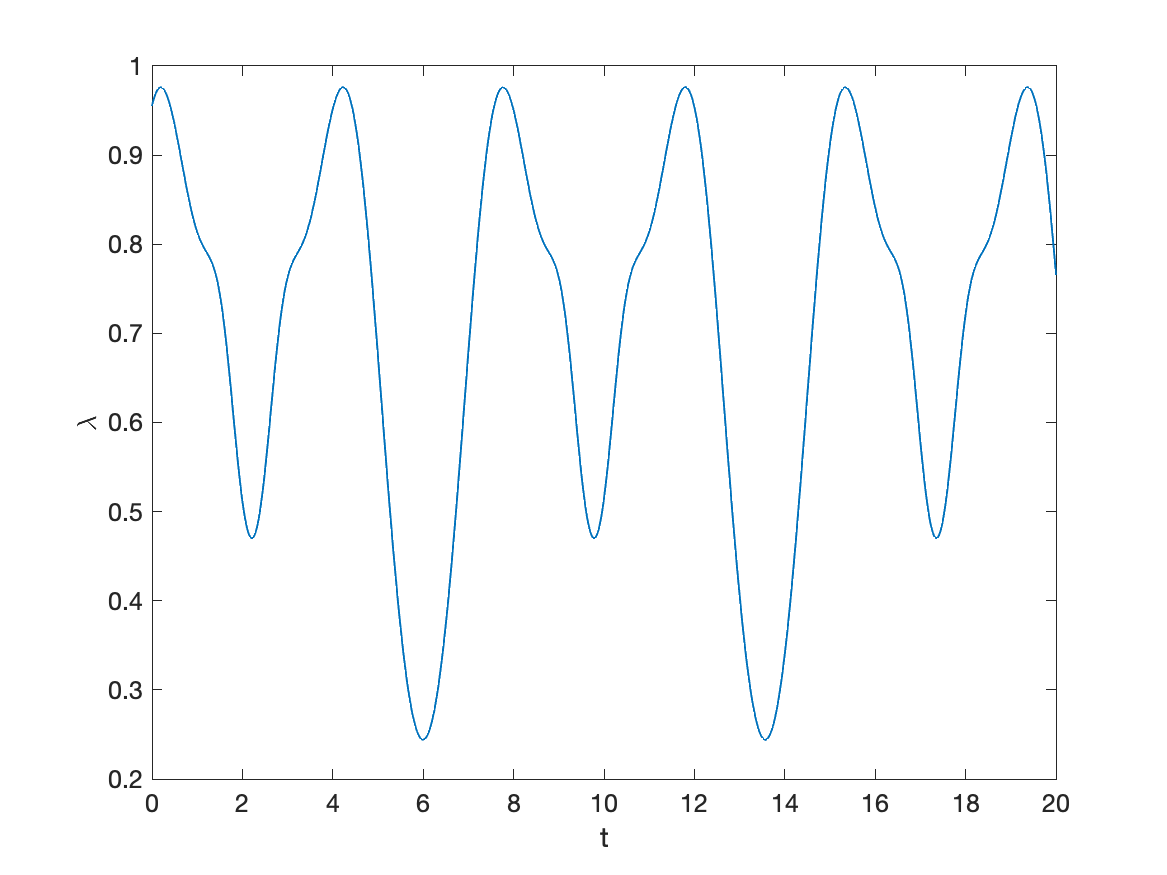}}
\caption{\label{pend_lam1} phase portrait (left) and Lagrange multiplier (right) for the problem defined by (\ref{pendHmod}), with $q^*=(2,0)^\top$, using the initial values (\ref{pend0}); in the left-plot, the star denotes the second charged point. }
\end{figure}

In this case, the simple pendulum described in Section~\ref{simpend} is modified by assuming that the unit-point mass also has an electric charge. Another charged point is placed at $q^*$, so that the Hamiltonian becomes:
\begin{equation}\label{pendHmod}
H(q,p) = \frac{1}2 \|p\|^2 +e_2^\top q -\frac{1}{\|q-q^*\|},
\end{equation}
subject to the same constraint in (\ref{pend}). We consider the same initial values (\ref{pend0}), and choose $q^*=(2,\,0)^\top$. The corresponding phase portrait and the Lagrange multiplier are depicted in Figure~\ref{pend_lam1}, relative to the interval $[0,20]$.

As is clear, the potential is now no more quadratic, so that the HBVM$(s,s)$ method is no more energy-conserving. Nevertheless, by using a HBVM$(k,s)$ method, with $k\ge s$, the Hamiltonian error should behave as $O(h^{2k})$, according to Theorem~\ref{evaik}. This is confirmed by the numerical results listed in Table~\ref{pendtab1}, obtained by solving the problem on the interval $[0,20]$, with timesteps as defined in (\ref{hi}),
by means of the HBVM$(k,1)$ method, $k=1,2,3,4$. As one may see, the HBVM(4,1) method is (practically) energy-conserving, for the considered values of the timestep.

\begin{table}[t]
\caption{obtained results for problem (\ref{pend0}) and (\ref{pendHmod}) by using HBVM$(k,1)$ methods with timestep (\ref{hi}). The *** mean that the round-off error level is reached.}\label{pendtab1}
\smallskip\small
\centerline{
\begin{tabular}{|cccccccccc|}
\hline
\multicolumn{10}{|c|}{$k=1$}\\
\hline
$i$ & $e_y$ & rate & $e_\lambda$ & rate & $e_{hid}$ & rate & $e_g$ & $e_H$ & rate\\
\hline
  3 & 2.40e-02 &   ---  & 3.78e-02 &    ---  & 1.25e-03 &   ---  & 2.00e-15 &    9.10e-05 &   --- \\ 
  4 & 5.90e-03 &   2.0 & 1.52e-02 &   1.3 & 3.12e-04 &   2.0 & 1.62e-14 &    2.28e-05 &   2.0 \\ 
  5 & 1.47e-03 &   2.0 & 6.54e-03 &   1.2 & 7.80e-05 &   2.0 & 1.78e-15 &   5.69e-06 &   2.0 \\ 
  6 & 3.67e-04 &   2.0 & 3.00e-03 &   1.1 & 1.95e-05 &   2.0 & 2.89e-15 &   1.42e-06 &   2.0 \\ 
  7 & 9.17e-05 &   2.0 & 1.43e-03 &   1.1 & 4.87e-06 &   2.0 & 5.77e-15 &    3.56e-07 &   2.0 \\ 
\hline
\hline
\multicolumn{10}{|c|}{$k=2$}\\
\hline
$i$ & $e_y$ & rate & $e_\lambda$ & rate & $e_{hid}$ & rate & $e_g$ & $e_H$ & rate\\
\hline
  3 & 2.28e-02 &   ---  & 3.70e-02 &   ---  & 1.25e-03 &   ---   & 2.66e-15 &    8.65e-09 &   --- \\ 
  4 & 5.63e-03 &   2.0 & 1.50e-02 &   1.3 & 3.12e-04 &   2.0 & 1.82e-14 &   5.41e-10 &   4.0 \\ 
  5 & 1.40e-03 &   2.0 & 6.49e-03 &   1.2 & 7.80e-05 &   2.0 & 4.44e-15 &    3.38e-11 &   4.0 \\ 
  6 & 3.50e-04 &   2.0 & 2.99e-03 &   1.1 & 1.95e-05 &   2.0 & 4.11e-15 &    2.11e-12 &   4.0 \\ 
  7 & 8.76e-05 &   2.0 & 1.43e-03 &   1.1 & 4.87e-06 &   2.0 & 3.33e-15 &    1.25e-13 &   4.1 \\ 
\hline
\hline
\multicolumn{10}{|c|}{$k=3$}\\
\hline
$i$ & $e_y$ & rate & $e_\lambda$ & rate & $e_{hid}$ & rate & $e_g$ & $e_H$ & rate\\
\hline
  3 & 2.28e-02 &   --- & 3.70e-02 &   --- & 1.25e-03 &   --- & 2.11e-15 &    7.32e-13 &   --- \\ 
  4 & 5.63e-03 &   2.0 & 1.50e-02 &   1.3 & 3.12e-04 &   2.0 & 5.55e-16 &   9.77e-15 &   6.2 \\ 
  5 & 1.40e-03 &   2.0 & 6.49e-03 &   1.2 & 7.80e-05 &   2.0 & 5.55e-16 &   5.55e-16 &   *** \\ 
  6 & 3.50e-04 &   2.0 & 2.99e-03 &   1.1 & 1.95e-05 &   2.0 & 7.77e-16 &   6.66e-16 &  *** \\ 
  7 & 8.76e-05 &   2.0 & 1.43e-03 &   1.1 & 4.87e-06 &   2.0 & 2.89e-15 &   1.11e-15 &  *** \\ 
\hline
\hline
\multicolumn{10}{|c|}{$k=4$}\\
\hline
$i$ & $e_y$ & rate & $e_\lambda$ & rate & $e_{hid}$ & rate & $e_g$ & $e_H$ & rate\\
\hline
  3 & 2.28e-02 &   --- & 3.70e-02 &   --- & 1.25e-03 &   --- & 2.22e-16 &       4.44e-16 &   --- \\ 
  4 & 5.63e-03 &   2.0 & 1.50e-02 &   1.3 & 3.12e-04 &   2.0 & 4.44e-16 &   4.44e-16 &  *** \\ 
  5 & 1.40e-03 &   2.0 & 6.49e-03 &   1.2 & 7.80e-05 &   2.0 & 4.44e-16 &   6.66e-16 &  *** \\ 
  6 & 3.50e-04 &   2.0 & 2.99e-03 &   1.1 & 1.95e-05 &   2.0 & 6.66e-16 &   6.66e-16 &  *** \\ 
  7 & 8.76e-05 &   2.0 & 1.43e-03 &   1.1 & 4.87e-06 &   2.0 & 1.55e-15 &   1.11e-15 &  *** \\ 
\hline
\hline
\end{tabular}}
\end{table}

\subsection{A conical pendulum}\label{conpend}

\begin{figure}[h]
\centering
\includegraphics[width=9cm]{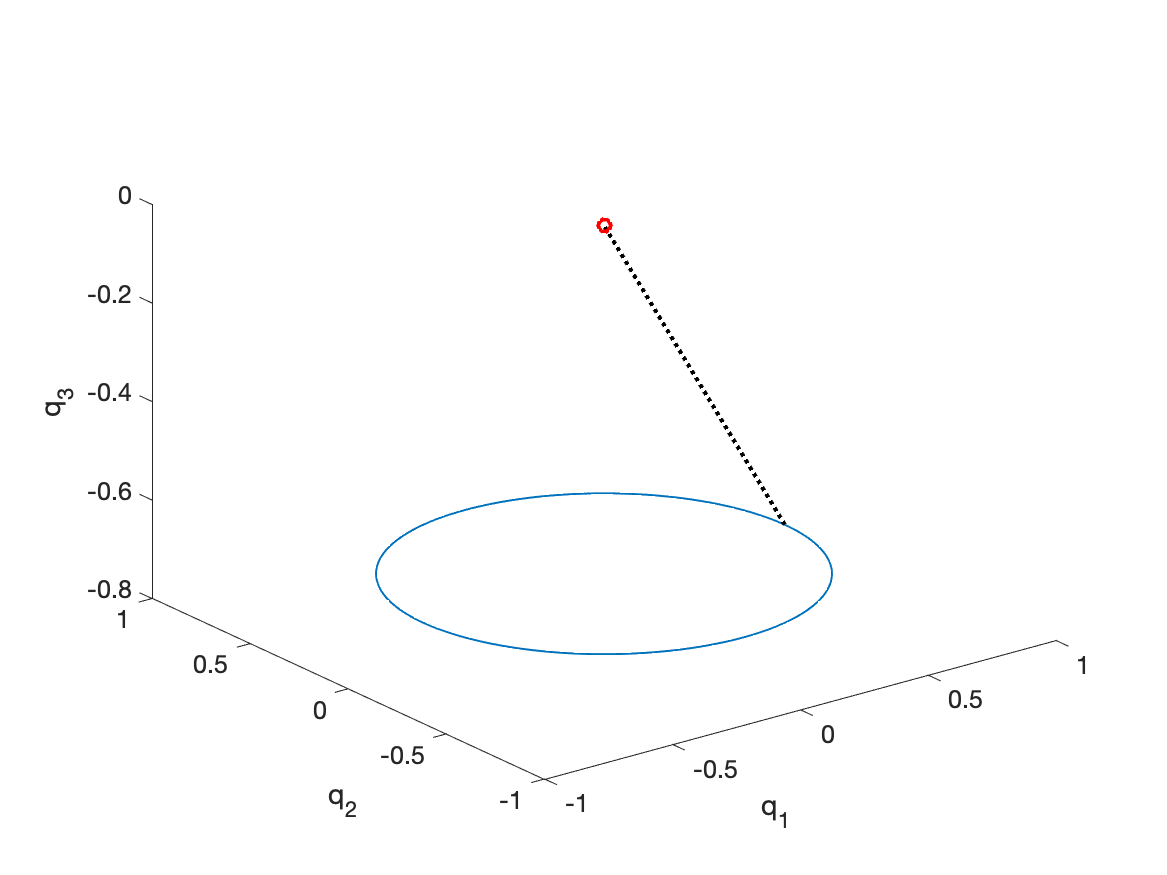}
\caption{\label{conpend_lam} phase portrait the problem  (\ref{cpend})--(\ref{cpend0}). }
\end{figure}

\begin{table}[t]
\caption{obtained results for problem  (\ref{cpend})--(\ref{cpend0}) by using HBVM$(s,s)$ methods with timestep (\ref{hi2}). The *** mean that the round-off error level is reached.}\label{cpendtab}
\smallskip\small
\centerline{
\begin{tabular}{|ccccccc|}
\hline
\multicolumn{7}{|c|}{$s=1$}\\
\hline
$i$ & $e_y$ & rate & $e_\lambda$  & $e_{hid}$  & $e_g$ & $e_H$ \\
\hline
  0 & 3.61e-01 &   --- & 3.03e-14 &   3.33e-16 &   7.55e-15 &   2.22e-15  \\ 
  1 & 1.20e-01 &   1.6 & 3.33e-16 &   1.28e-15 &   5.55e-16 &   1.72e-15  \\ 
  2 & 3.20e-02 &   1.9 & 1.67e-15 &   4.16e-16 &   7.77e-16 &   2.22e-16  \\ 
  3 & 8.11e-03 &   2.0 & 1.60e-14 &   2.33e-15 &   4.44e-16 &   1.67e-16  \\ 
  4 & 2.03e-03 &   2.0 & 1.57e-14 &   1.33e-15 &   7.77e-16 &   4.44e-16  \\ 
\hline
\hline
\multicolumn{7}{|c|}{$s=2$}\\
\hline
$i$ & $e_y$ & rate & $e_\lambda$  & $e_{hid}$  & $e_g$ & $e_H$ \\
\hline
  0 & 1.55e-02 &   --- & 1.51e-13 &   1.11e-15 &   1.11e-15 &   6.66e-16  \\ 
  1 & 1.05e-03 &   3.9 & 2.66e-15 &  1.33e-15 &   1.33e-15 &  1.39e-15  \\ 
  2 & 6.66e-05 &   4.0 & 1.07e-14 &  1.69e-15 &   2.22e-16 &   2.78e-16  \\ 
  3 & 4.18e-06 &   4.0 & 1.25e-14 &  1.44e-15 &   4.44e-16 &   3.89e-16  \\ 
  4 & 2.62e-07 &   4.0 & 1.59e-13 &  4.39e-15 &   6.66e-16 &   3.33e-16  \\ 
\hline
\hline
\multicolumn{7}{|c|}{$s=3$}\\
\hline
$i$ & $e_y$ & rate & $e_\lambda$  & $e_{hid}$  & $e_g$ & $e_H$ \\
\hline
  0 & 1.91e-04 &   ---  & 2.00e-15 &  9.44e-16 &  3.33e-16 &  1.11e-16  \\ 
  1 & 3.13e-06 &   5.9 & 6.22e-15 &  9.99e-16 &  3.33e-16 &  1.11e-16  \\ 
  2 & 4.94e-08 &   6.0 & 1.18e-14 &  6.66e-16 &  5.55e-16 &  2.22e-16  \\ 
  3 & 7.74e-10 &   6.0 & 2.24e-14 &  8.88e-16 &  8.88e-16 &  2.22e-16  \\ 
  4 & 1.21e-11 &   6.0 & 2.24e-13 &  2.44e-15 &  8.88e-16 &  6.66e-16  \\ 
\hline
\hline
\multicolumn{7}{|c|}{$s=4$}\\
\hline
$i$ & $e_y$ & rate & $e_\lambda$  & $e_{hid}$  & $e_g$ & $e_H$ \\
\hline
  0 & 1.23e-06 &    --- & 1.02e-14 &   2.75e-15 &   2.22e-16 &   4.44e-16  \\ 
  1 & 4.97e-09 &   7.9 & 1.98e-14 &   6.66e-16 &   6.66e-16 &   5.55e-16  \\ 
  2 & 1.96e-11 &   8.0 & 9.18e-14 &   2.56e-15 &   4.44e-16 &   1.11e-16  \\ 
  3 & 7.69e-14 &   8.0 & 1.38e-13 &   3.08e-15 &   4.44e-16 &   2.78e-16 \\ 
  4 & 8.80e-16 &    *** & 6.98e-14 &   2.00e-15 &   4.44e-16 &   3.89e-16  \\ 
\hline
\hline
\end{tabular}}
\end{table}

In this case, the problem is given by:
\begin{equation}\label{cpend}
\dot q = p, \qquad \dot p = -e_3, \qquad g(q) := \|q\|^2-1=0, 
\end{equation}
$e_3$ being the third unit vector in $\RR^3$, with Hamiltonian
\begin{equation}\label{cpendH}
H(q,p) = \frac{1}2 \|p\|^2 +e_3^\top q.
\end{equation}
Considering the set of consistent initial conditions:
\begin{equation}\label{cpend0}
q(0) = \pmatrix{c} z_0\\  0 \\ -z_0\endpmatrix, \qquad p(0) = \pmatrix{c}0\\ \sqrt{z_0}\\ 0\endpmatrix, \qquad z_0 = \frac{1}{\sqrt{2}},
\end{equation}
the motion lies in the plane $q_3=-z_0$, with period $T=2^{\frac{3}4}\pi$, and the Lagrange multiplier is constant, $\lambda = z_0$. The suspension rod then spans a cone, from which the name of the problem (see Figure~\ref{conpend_lam}). For this problem, the Lagrange multiplier is exactly computed, so that condition (\ref{dli0}) holds true and HBVM$(s,s)$ methods are expected to be energy-conserving and with order $2s$. This is confirmed by the numerical tests listed in Table~\ref{cpendtab}, obtained by using the timestep
\begin{equation}\label{hi2}
h_i = 2^{-i}\frac{T}5, \qquad i=0,1,\dots,
\end{equation}
for solving one period. From the results in Table~\ref{cpendtab}, one may see that remarkably enough, also the hidden constraint turns out to be satisfied.


\subsection{A tethered satellite system}\label{tetsys}

\begin{figure}[p]
\centering
\includegraphics[width=6.5cm]{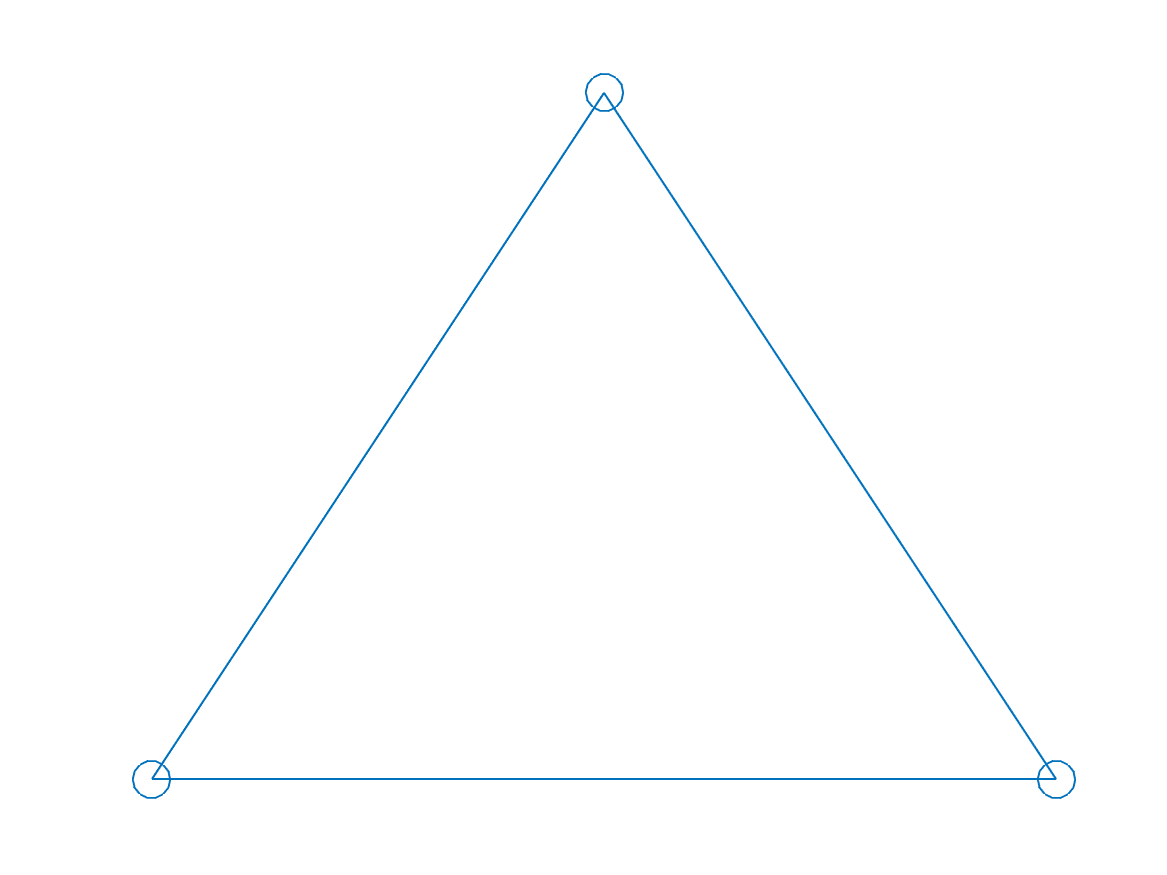} 
\caption{\label{tetfig0} tethered satellite system described by problem (\ref{tetH})--(\ref{p0}).}

\bigskip 
\centerline{
\includegraphics[width=6.5cm]{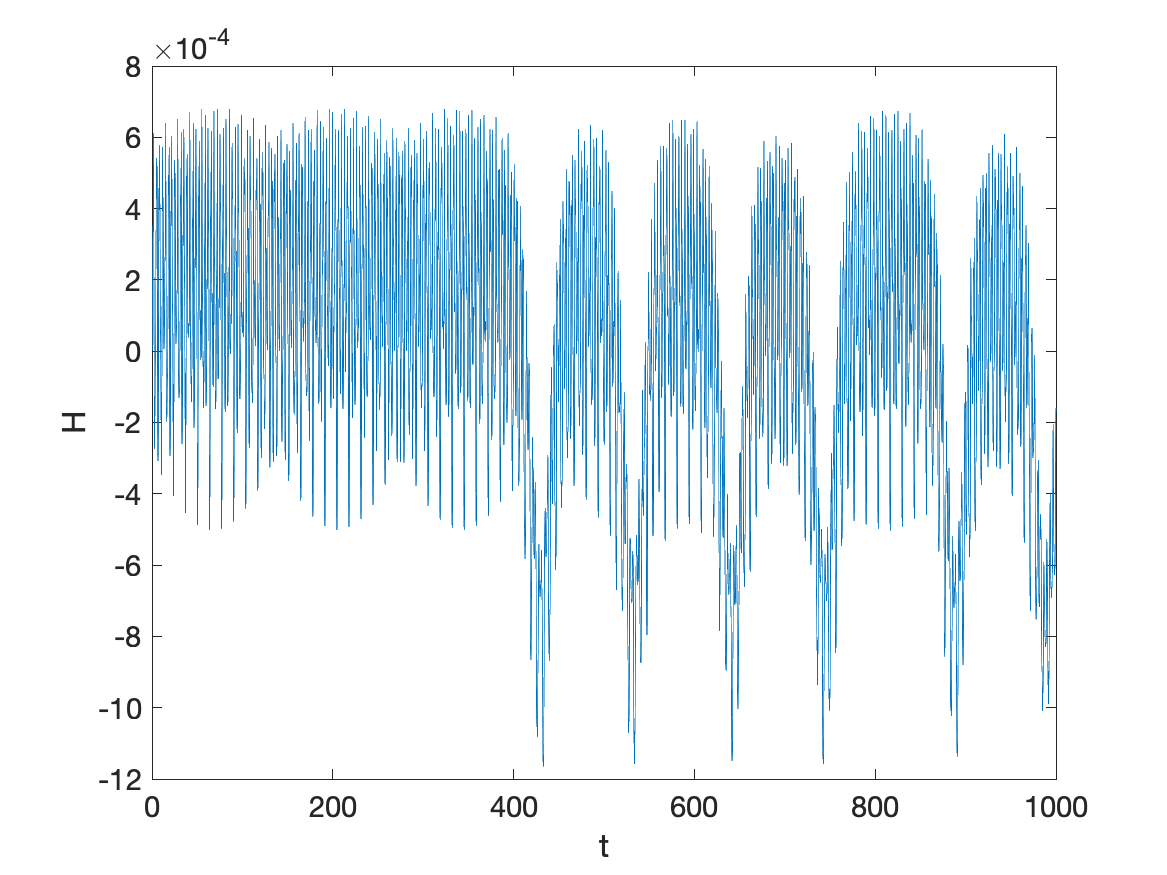}\quad \includegraphics[width=6.5cm]{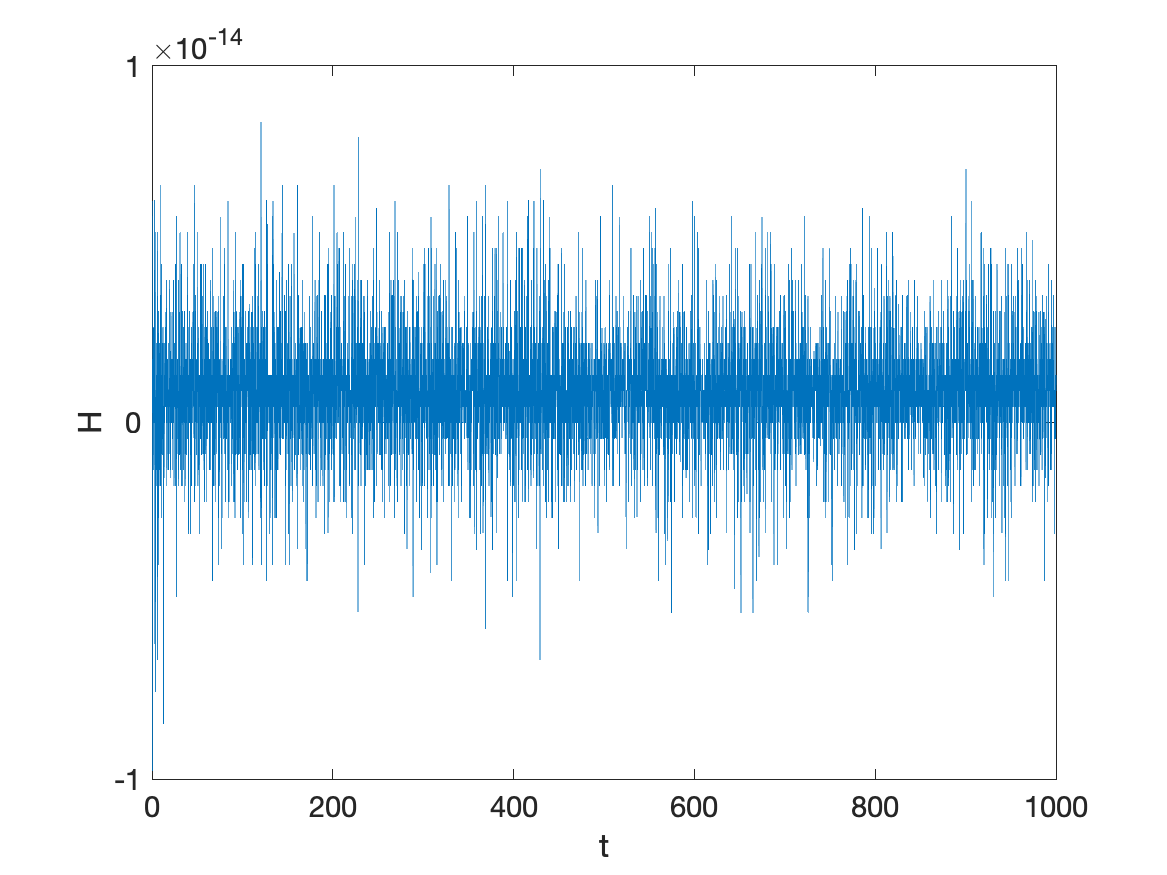}}
\caption{\label{tetfig1} numerical Hamiltonian when solving problem (\ref{tetH})--(\ref{p0}) with timestep $h=0.1$ by using HBVM(1,1) (left plot) and HBVM(5,1) (right plot). }

\bigskip
\centerline{
\includegraphics[width=6.5cm]{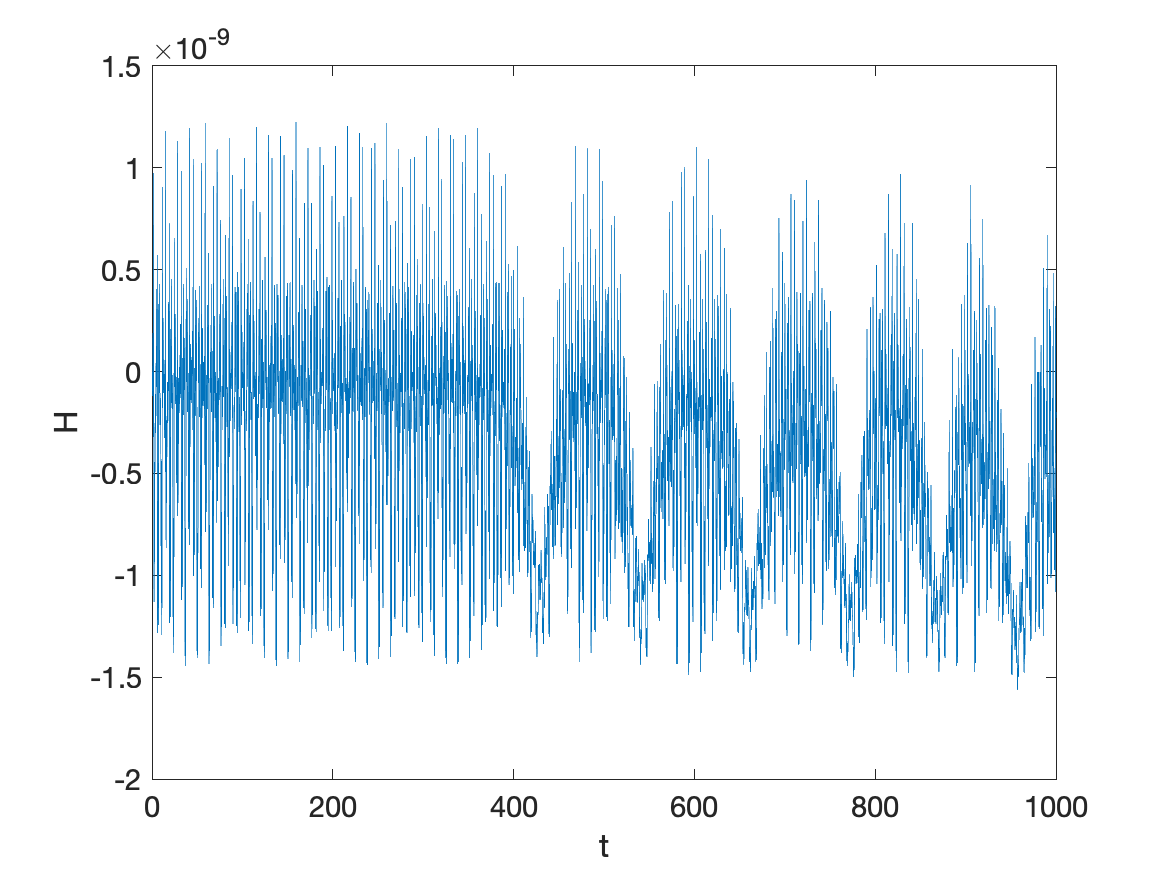}\quad \includegraphics[width=6.5cm]{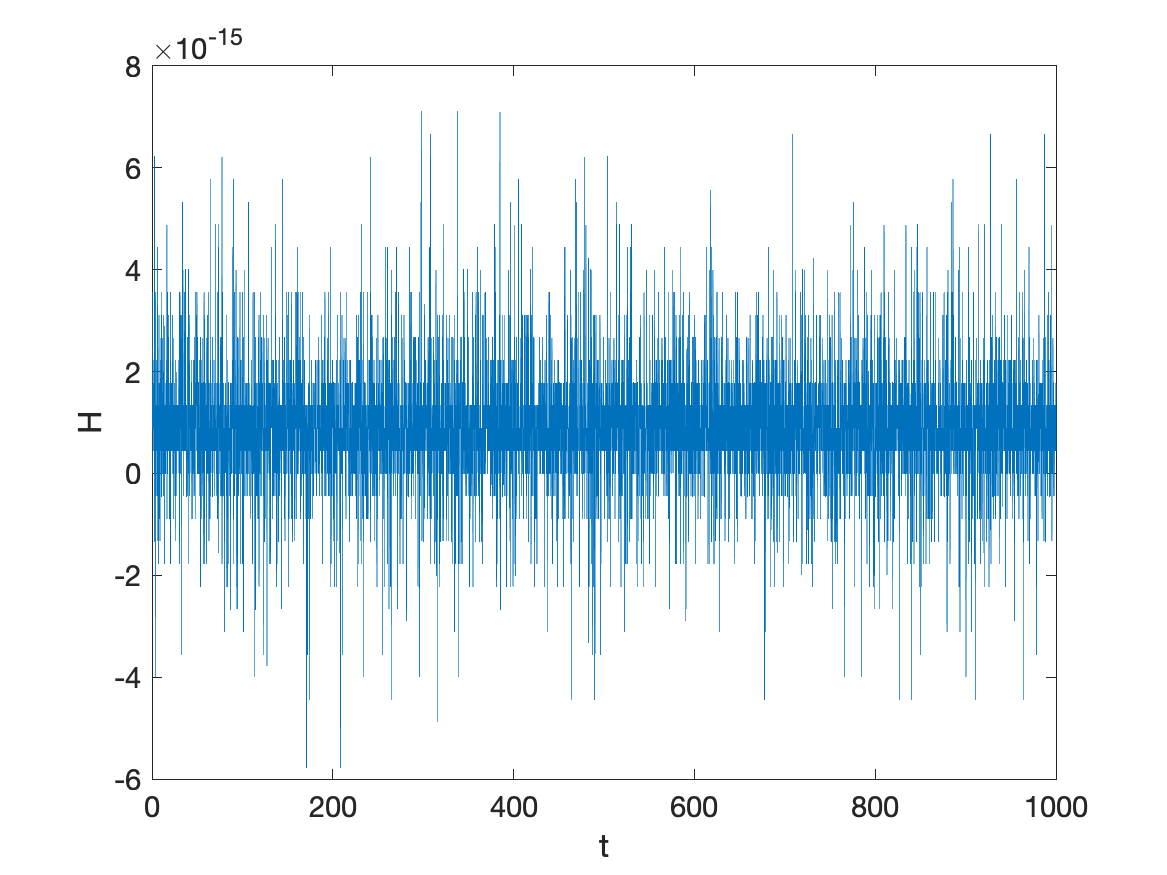}}
\caption{\label{tetfig2} numerical Hamiltonian when solving problem (\ref{tetH})--(\ref{p0}) with timestep $h=0.1$ by using HBVM(3,3) (left plot) and HBVM(5,3) (right plot). }

\end{figure}

The last test problem that we consider is a tethered satellite systems of three satellites, forming an equilateral triangle with unit sides, as depicted in Figure~\ref{tetfig0}. Setting $$q = \pmatrix{c}q_1\\ q_2\\ q_3\endpmatrix, \qquad p = \pmatrix{c}p_1\\ p_2\\ p_3\endpmatrix, \qquad q_i,p_i\in\RR^3, \quad i=1,2,3,$$
we assume the Hamiltonian be given by
\begin{equation}\label{tetH}
H(q,p) = \sum_{i=1}^3 \left(\frac{1}2 \|p_i\|^2 +\frac{1}{\|q_i\|} +\cos(\|q_i\|)\right),
\end{equation}
with the constraints:
\begin{equation}\label{gi}
g(q) := \pmatrix{c}\|q_1-q_2\|^2-1 \\ \|q_2-q_3\|^2-1 \\ \|q_1-q_3\|^2-1\endpmatrix = 0\in\RR^3.
\end{equation}
Consistent initial conditions are chosen as follows:
\begin{equation}\label{q0}
q_1(0) = \pmatrix{c} 0\\ \frac{1}2 \\ z_0 \endpmatrix, \quad
q_2(0) = \pmatrix{c} 0\\ -\frac{1}2 \\ z_0 \endpmatrix, \quad
q_3(0) = \pmatrix{c} 0\\ 0 \\ z_0-\frac{\sqrt{3}}2 \endpmatrix, 
\end{equation}
and
\begin{equation}\label{p0}
p_1(0) = p_2(0) = 0\in\RR^3, \qquad p_3(0) = \pmatrix{c} v_0 \\ 0\\ 0\endpmatrix,
\end{equation}
with $z_0=20$ and $v_0>0$ such that $H(q(0),p(0))=0$.
Let us numerically solve problem (\ref{tetH})--(\ref{p0}) on the interval $[0,10^3]$ by using the HBVM(1,1) method with timestep $h=0.1$. The method is not energy-conserving, and 9.4 sec are required for solving the problem, with a mean number of 15.4 iterations (\ref{iter}) per step. The numerical Hamiltonian is depicted in the left-plot of Figure~\ref{tetfig1}. If we solve the problem by using HBVM(5,1), the method turns out to be practically energy-conserving, and requires 11.6 sec, with the same mean number of iterations (\ref{iter}) per step. The corresponding numerical Hamiltonian is shown in the right-plot of Figure~\ref{tetfig1}. As one may see, the energy-conserving method has a much more appropriate behavior, with a very moderate increase of the computational time.

Similar results are obtained by solving the problem, with the same timestep $h=0.1$, by using the HBVM(3,3) and HBVM(5,3) methods, the latter method being (practically) energy-conserving. Both methods need a mean of 14 iterations (\ref{iter}) per step, with execution times of 14.2 sec and 15.4 sec, respectively. The corresponding numerical Hamiltonians are depicted in Figure~\ref{tetfig2}.  Also in this case, the qualitative results obtained by the energy-conserving method are more favorable, with only a very limited increase of the overall computational cost.

\section{Conclusions}\label{fine}

In this paper we have introduced a class of energy-conserving methods for Hamiltonian problems with quadratic holonomic constraints. The methods reduce to the class of Hamiltonian Boundary Value Methods (HBVMs), when the constraints are not given. Arbitrarily high-order methods are obtained within this class. A detailed analysis of the methods has been presented, along with their practical implementation and solution of the generated discrete problems. A number of numerical tests duly confirm the theoretical results, also proving the usefulness of the approach. 
Future directions of investigation will include the extension of the methods for numerically solving Hamiltonian problems with general holonomic constraints.

\paragraph*{Conflict of interests.} The authors declare no conflict of interests.

\paragraph*{Acknowledgements.} The authors acknowledge the {\em mrSIR project} \cite{mrsir} for the financial support.

\end{document}